\documentclass[12pt]{article}
\usepackage{amssymb}
\usepackage{latexsym}
\def\part#1{\frac{\partial\phantom{q}}{\partial#1}}

\newenvironment{rmk}{\begin{trivlist}\item[]{\bf Remark:} }
{\end{trivlist}}
\newenvironment{ex}{\begin{trivlist}\item[]{\bf Example:} }
{\end{trivlist}}
\newenvironment{prf}{\begin{trivlist}\item[]{\bf Proof:} }
{\hfill $\Box$ \end{trivlist}}
\newenvironment{lemprf}{\begin{trivlist}\item[]{\bf Proof:} }
 {\end{trivlist}}
\newtheorem{thm}{Theorem}
\newtheorem{definition}{Definition}
\newtheorem{prp}[thm]{Proposition}
\newtheorem{lem}[thm]{Lemma}
\newtheorem{cor}[thm]{Corollary}

\newcommand{\lie}[1]{\mathfrak{#1}}
\def\End{\mathop{\rm End}\nolimits}

\def\ker{\mathop{\rm ker}\nolimits}

\def\deg{\mathop{\rm deg}\nolimits}
\def\rk{\mathop{\rm rk}\nolimits}
\def\tr{\mathop{\rm tr}\nolimits}

\def\star{\mathop{*\!}\nolimits}
\def\cyc{\mathrm {cyclic}}
\def\dstar{\mathop{d_A*\!}\nolimits}
\newcommand{\R}{\mathbf{R}}
\newcommand{\C}{\mathbf{C}}

\newcommand{\CP}{{\mathbf C}{\rm P}}
\newcommand{\RP}{{\mathbf R}{\rm P}}
\begin{document}
\title{Instantons, Poisson structures  and generalized K\"ahler geometry}
 \author{Nigel Hitchin\\[5pt]
\itshape  Mathematical Institute\\
\itshape 24-29 St Giles\\
\itshape Oxford OX1 3LB\\
\itshape UK\\
 hitchin@maths.ox.ac.uk}
\maketitle
\begin{abstract}
\noindent Using the idea of a generalized K\"ahler structure, we construct bihermitian metrics on $\CP^2$ and $\CP^1\times \CP^1$, and show that  any such structure on a compact $4$-manifold $M$ defines one on the moduli space of anti-self-dual connections on a fixed principal bundle over $M$. We highlight the role of holomorphic Poisson structures in all these constructions.
 \end{abstract}
\section{Introduction}
The idea of a generalized complex structure -- a concept which interpolates between complex and symplectic structures -- seems to provide a differential geometric language in which some of the structures of current interest in string theory fit very naturally. There is an associated notion of \emph{generalized K\"ahler manifold} which essentially consists of a pair of commuting generalized complex structures. A remarkable theorem of Gualtieri \cite{Gu} shows that it has an equivalent interpretation in standard geometric terms: a manifold with two complex structures $I_+$ and $I_-$; a metric $g$, Hermitian with respect to both; and connections $\nabla_+$ and $\nabla_-$ compatible with these structures but with skew torsion $db$ and $-db$ respectively for a $2$-form $b$. This so-called \emph{bihermitian structure} appeared in the physics literature as long ago as 1984 \cite{R} as a target space for the supersymmetric $\sigma$-model and in the pure mathematics literature more recently (\cite{AGG} for example) in the context of the integrability of the canonical almost complex structures defined by the Weyl tensor of a Riemannian four-manifold. The theory has suffered from a lack of interesting examples.

The first purpose of this paper is to use the generalized complex structure approach to find non-trivial explicit  examples on $\CP^2$ and $\CP^1\times \CP^1$. We use an approach to generalized K\"ahler structures of generic type which involves closed $2$-forms satisfying algebraic conditions. This is in principle much easier than trying to write down the differential-geometric data above. What we show is that every $SU(2)$-invariant K\"ahler metric on $\CP^2$ or the Hirzebruch surface ${\mathbf F}_2$ generates naturally a generalized K\"ahler structure, where for ${\mathbf F}_2$ (which is diffeomorphic to $S^2\times S^2$), the complex structures $I_+,I_-$ are equivalent to ${\mathbf F}_0=\CP^1\times \CP^1$.  

The second part of the paper shows that a bihermitian structure on a $4$-manifold (where $I_+$ and $I_-$ define the same orientation) defines naturally a bihermitian structure  on the moduli space of solutions to the anti-self-dual Yang-Mills equations, and this gives another  (less explicit) source of examples.

What appears naturally in approaching these goals is the appearance of holomorphic Poisson structures, and in a way the main point of the paper is to bring this aspect into the foreground. It seems as if this type of differential geometry is related to complex Poisson manifolds in the way in which hyperk\"ahler metrics are adapted to complex symplectic manifolds. Yet our structures are more flexible -- like K\"ahler metrics they can be changed in the neighbourhood of a point. The link with Poisson geometry occurs in three   interlinking ways:
\begin{itemize}
\item
a holomorphic Poisson structure defines a particular type of generalized complex structure (see \cite{Gu}),
\item
the skew form $g([I_+,I_-]X,Y)$ for the bihermitian metric is of type $(2,0)+(0,2)$ and defines a holomorphic Poisson structure for either complex structure $I_+$ or $I_-$ (in the four-dimensional case this was done in \cite{AGG}),
\item
a generalized complex structure $J:T\oplus T^*\rightarrow T\oplus T^*$ defines by restriction a homomorphism $\pi:T^*\rightarrow T$ which is a real Poisson structure (this has been noted by several authors, see \cite{AB}).

\end{itemize}
We may also remark that Gualtieri's deformation theorem \cite{Gu} showed that interesting deformations of complex manifolds as generalized complex manifolds require the existence of a holomorphic Poisson structure.

We address all three Poisson-related issues in the paper. The starting point for our examples is the generalized complex structure determined by a complex Poisson surface (namely, a surface with an anticanonical divisor) and we solve the equations for a second generalized complex structure which commutes with this one. When we study the moduli space of instantons we show that the holomorphic Poisson structures defined by $g([I_+,I_-]X,Y)$ are the canonical ones studied by Bottacin \cite{Bot1}. Finally we examine the symplectic leaves of the real Poisson structures $\pi_1,\pi_2$ on the moduli space.

The structure of the paper is as follows. We begin by studying generalized K\"ahler manifolds as a pair $J_1,J_2$ of commuting generalized complex structures, and we focus in particular on the case where each $J_1,J_2$ is the B-field transform of a symplectic structure -- determined by a closed form $\exp (B+i\omega)$ --  giving a convenient algebraic form for the commuting property. We then implement this to find the two examples. In the next section we introduce the bihermitian interpretation and prove that $g([I_+,I_-]X,Y)$ does actually define a holomorphic Poisson structure. 

The following sections show how to introduce a bihermitian structure on the moduli space ${\mathcal M}$ of gauge-equivalence classes of solutions to the anti-self-dual Yang-Mills equations. At first glance this seems obvious -- we have two complex structures $I_+,I_-$ on $M$ and hence two complex structures on ${\mathcal M}$, since ${\mathcal M}$ is the moduli space of $I_+$- or $I_-$-stable bundles, and we have a natural ${\mathcal L}^2$ metric. This would be fine for a K\"ahler metric but not in the non-K\"ahler case. Here L\"ubke and Teleman \cite{LT} reveal the correct approach -- one chooses a  different horizontal to the gauge orbits in order to define the metric on the quotient. In our case we have two complex structures and  two horizontals and much of the manipulation and integration by parts which occurs in this paper is caused by this complication. 

One aspect we do not get is a natural pair of commuting generalized complex structures on ${\mathcal M}$ -- we obtain the differential geometric data above, and an exact $3$-form $db$, but not a natural choice of $b$.  We get a generalized K\"ahler structure only modulo a closed B-field on ${\mathcal M}$. This suggests that ${\mathcal M}$ is not, at least directly, a moduli space of objects  defined solely by one of the commuting generalized complex structures on $M$, but there is clearly more to do here. 

We give finally a quotient construction  which also demonstrates the problem of making a generalized K\"ahler structure descend to the quotient. This procedure, analogous to the hyperk\"ahler quotient,  could be adapted to yield the bihermitian metric on ${\mathcal M}$ in the case $M$ is a $K3$ or torus.  Unfortunately we have not found a quotient construction for the instanton moduli space which works in full generality, but this might be possible by using framings on the anticanonical divisor.

\vskip .5cm
\noindent{{\bf Acknowledgements:}} The author wishes to thank M. Gualtieri, G. Cavalcanti and V. Apostolov for useful discussions.

\section{Generalized K\"ahler manifolds}
\subsection{Basic properties}

The notion of a generalized K\"ahler structure was introduced by M. Gualtieri in \cite{Gu}, in the context of the generalized complex structures defined by the author in \cite{Hit}. Recall that ``generalized geometry" consists essentially of replacing the tangent bundle $T$ of a manifold by $T\oplus T^*$ with its natural indefinite inner product $$(X+\xi,X+\xi)=-i_X\xi,$$ and  the Lie bracket on sections of $T$ by the Courant bracket 
$$[X+\xi, Y+\eta]=
  [X,Y]+\mathcal{L}_{X}\eta -\mathcal{L}_{Y}\xi -\frac{1}{2}
  d(i_{X}\eta  -i_{Y}\xi)$$
  on sections of $T\oplus T^*$. One then  introduces additional structures on $T\oplus T^*$ compatible with these. 
A \emph{generalized complex structure} is a complex structure $J$ on $T\oplus T^*$ such that $J$ is orthogonal with respect to the inner product and with the integrability condition that if $A,B$ are sections of $(T\oplus T^*)\otimes \C$ with $JA=iA,JB=iB$, then $J[A,B]=i[A,B]$ (using the Courant bracket). The standard examples are a complex manifold 
where $$J_1=\pmatrix {I&0\cr
                                0& -I}$$
                               and a 
                                symplectic manifold where $$J_2=\pmatrix {0&-\omega^{-1}\cr
                                \omega & 0}.$$
                                
 The  $+i$ eigenspace  of $J$   is spanned by     $\{\dots, \partial/\partial z_j\dots, \dots, d\bar z_k,\dots\}$  in the first case  and $\{\dots, \partial/\partial x_j- i \sum\omega_{jk}dx_k,\dots\}$ in the second. 

Another example of a generalized complex manifold is a holomorphic Poisson manifold -- a complex manifold with a holomorphic bivector field
$$\sigma=\sum\sigma^{ij}\frac{\partial}{\partial  z_i}\wedge\frac{\partial}{\partial  z_j}$$
satisfying the condition $[\sigma,\sigma]=0$, using the Schouten bracket. This defines a generalized complex structure where the $+i$ eigenspace is 
$$ E=\left[\dots,  \frac{\partial}{\partial  z_j},\dots, d\bar z_k+\sum_{\ell}\bar\sigma^{k\ell}\frac{\partial}{\partial  \bar z_{\ell}},\dots \right],$$    
and if $\sigma=0$ this gives a complex structure.                   
                              
     Gualtieri   observed  (see also \cite{AB}) that the real bivector defined by the upper triangular part of  $J:T\oplus T^*\rightarrow T\oplus T^*$ is always a real Poisson structure. In the symplectic case this is the canonical Poisson structure and in the complex case it is zero. Both facts show that  Poisson geometry  plays a central role in this area, a feature we shall see more of later.
\vskip .25cm                       
   The algebraic compatibility condition between $\omega$ and $I$ to give a K\"ahler manifold (i.e. that $\omega$ be of type $(1,1)$) can be expressed as $J_1J_2=J_2J_1$   and this is the basis of the definition of a \emph{generalized}  K\"ahler structure:

\begin{definition} \label{GKdef} A   \emph{generalized K\"ahler structure} on a manifold consists of two commuting generalized complex structures  $J_1,J_2$ such that the quadratic form $(J_1J_2A,A)$ on $T\oplus T^*$ is  definite.  
\end{definition}            
\vskip .25cm
At a point, a generalized complex structure can also be described  by  a form $\rho$: the $+i$ eigenspace bundle $E$  consists of the $A=X+\xi\in (T\oplus T^*)\otimes \C$ which satisfy $A\cdot \rho=i_X\rho+\xi\wedge\rho=0$. For the symplectic structure $\rho=\exp i\omega$, and for a complex structure with complex coordinates $z_1,\dots,z_n$ we take the $n$-form $\rho=dz_1\wedge dz_2\wedge \dots\wedge dz_n$. The structure is called even or odd according to whether $\rho$ is an even or odd form. The generic even case is the so-called B-field transform of a symplectic structure where
$$\rho=\exp \beta = \exp (B+i\omega)$$
and $B$ is an arbitrary  $2$-form. The generalized complex structure defined by a holomorphic Poisson structure $\sigma$ is of this type if $\sigma$ defines a non-degenerate skew form on $(T^*)^{1,0}$; then $B+i\omega$ is its inverse.

If $\rho$ extends smoothly to a neighbourhood of the point, and is {\it closed}, then the integrability condition for a generalized complex structure holds.
\vskip .25cm
The following lemma is  useful for finding generalized K\"ahler structures where both are of this generic even type (which requires the dimension of $M$ to be of the form $4k$). We shall return to this case periodically to see how the various structures emerge concretely.

\begin{lem} \label{commute} Let  $\rho_1=\exp \beta_1, \rho_2=\exp\beta_2$ be closed forms defining generalized complex structures $J_1,J_2$ on a manifold of dimension $4k$. Suppose that
$$(\beta_1-\beta_2)^{k+1}=0=(\beta_1-\bar\beta_2)^{k+1}$$
and $(\beta_1-\beta_2)^{k}$ and $(\beta_1-\bar\beta_2)^{k}$ are non-vanishing. Then  $J_1$ and $J_2$ commute.
\end{lem}

\begin{lemprf} Suppose that $(\beta_1-\beta_2)^{k+1}=0$ and  $(\beta_1-\beta_2)^{k}$ is non-zero. Then the $2$-form $\beta_1-\beta_2$ has rank $2k$, i.e. the dimension of the space of vectors $X$ satisfying $i_X(\beta_1-\beta_2)=0$ is $2k$. Since $i_X 1+\xi\wedge 1=0$ if and only if $\xi=0$, this means that the space of solutions $A=X+\xi$ to 
$$A\cdot \exp (\beta_1-\beta_2)=0=A\cdot 1$$
is $2k$-dimensional. Applying the invertible map $\exp \beta_2$, the same is true of solutions to 
$$A\cdot \exp \beta_1=0=A\cdot \exp \beta_2.$$
This is the intersection $E_1\cap E_2$ of the two $+i$ eigenspaces. Repeating for $\beta_1-\bar\beta_2$ we get $E_1\cap \bar E_2$ to be $2k$-dimensional. These two bundles are common eigenspaces of $(J_1, J_2)$ corresponding to the eigenvalues $(i,i)$ and $(i,-i)$ respectively. Together with their conjugates they decompose $(T\oplus T^*)\otimes \C$ into a direct sum of common eigenspaces of $J_1,J_2$, thus $J_1J_2=J_2J_1$ on every element.
\end{lemprf}

We also need to address the  definiteness of $(J_1J_2 A,A)$ in Definition 1. Let $V_+$ be the $-1$ eigenspace  of $J_1J_2$ (the notation signifies $J_1=+J_2$ on $V_+$). This  is
$$E_1\cap E_2\oplus \bar E_1\cap \bar E_2.$$
 If $X$ is a vector in the $2k$-dimensional space defined by  $i_X(\beta_1-\beta_2)=0$ then $A=X-i_X\beta_2$ satisfies $A\cdot \exp \beta_1=0=A\cdot \exp \beta_2$, i.e. $A\in E_1\cap E_2$. But then 
\begin{equation}
(A+\bar A,A+\bar A)=i_X\beta_2(\bar X)+i_{\bar X}\bar\beta_2(X)=(\beta_2-\bar\beta_2)(X,\bar X)
\label{posit}
\end{equation}
so we need to have this form to be  definite. Note that interchanging the roles of $\beta_1,\beta_2$, this is the same as $(\beta_1-\bar\beta_1)(X,\bar X)$ being definite.

\subsection{Hyperk\"ahler examples}

A hyperk\"ahler manifold $M$ of dimension $4k$ provides a simple example of a generalized K\"ahler manifold. Let $\omega_1,\omega_2,\omega_3$ be the three K\"ahler forms corresponding to the complex structures $I,J,K$ and set 
$$\beta_1=\omega_1+\frac{i}{2}(\omega_2-\omega_3),\quad \beta_2=\frac{i}{2}(\omega_2+\omega_3).$$
Then $\beta_1-\beta_2=\omega_1-i\omega_3$ is a $J$-holomorphic symplectic $2$-form and so clearly satisfies the conditions of Lemma \ref{commute}. Similarly $\beta_1-\bar\beta_2= \omega_1+i\omega_2$  is holomorphic symplectic for $K$. The vectors $X$ satisfying $i_X(\beta_1-\beta_2)=0$ are the $(0,1)$ vectors for $J$, and $\beta_2-\bar\beta_2=i(\omega_2+\omega_3)$ whose $(1,1)$ part with respect to $J$ is $i\omega_2$. Thus
$$(\beta_2-\bar\beta_2)(X,\bar X)=i\omega_2(X,\bar X)$$
which is positive definite. Thus a hyperk\"ahler manifold satisfies all the conditions to be generalized K\"ahler.
\vskip .25cm
 D. Joyce observed (see \cite{AGG}) that one can deform this example. Let $f$ be a smooth real function on $M$, and use the symplectic form $\omega_1$ to define a Hamiltonian vector field. Now integrate it to a one-parameter group of symplectic diffeomorphisms $F_t:M\rightarrow M$, so that $F_t^*\omega_1=\omega_1$. Define
$$\beta_1=\omega_1+\frac{i}{2}(\omega_2-F_t^*\omega_3),\quad \beta_2=\frac{i}{2}(\omega_2+F_t^*\omega_3),$$
and then
$$\beta_1-\beta_2=\omega_1-iF^*_t\omega_3=F^*_t(\omega_1-i\omega_3).$$
This is just the pull-back by a diffeomorphism of $\omega_1-i\omega_3$ so also satisfies the constraint of Lemma \ref{commute}. We also have
$\beta_1-\bar\beta_2= \omega_1+i\omega_2$ which is just the same as the hyperk\"ahler case, so both constraints hold. If $t$ is sufficiently small this will still give a positive definite metric.

This simple example at least shows the flexibility of the concept -- we can find a new structure from an arbitrary real function, somewhat analogous to the addition of $\partial\bar\partial f$ to a K\"ahler form. In the compact four-dimensional situation this type of structure restricts us to tori and K3 surfaces. We give next an explicit example on the projective plane.

\subsection {Example: the projective plane}  \label{cp2}

The standard $SU(2)$ action on $\C^2$ extends to $\CP^2$ and the invariant $2$-form $dz_1\wedge dz_2$ extends to a meromorphic form with a triple pole on  the line at infinity. Its inverse $\partial/\partial z_1\wedge \partial/\partial z_2$ is a holomorphic Poisson structure with a triple zero on the line at infinity. We shall take the generalized complex structure $J_1$ to be defined by this, and seek an $SU(2)$-invariant generalized complex structure $J_2$ defined by $\exp (B+i\omega)$ in such a way that the pair define a generalized K\"ahler structure. On $\C^2$ the Poisson structure is non-degenerate, so the generalized complex structure on that open set is defined by the closed form $\rho_1=\exp dz_1dz_2$.
\vskip .25cm
We begin by parametrizing $\C^2\setminus \{0
\}$ by $\R^+ \times SU(2)$:
$$\pmatrix {z_1\cr
            z_2}=\pmatrix {z_1& -\bar z_2\cr
            z_2 & \bar z_1}\pmatrix {1\cr
            0}=rA\pmatrix {1\cr
            0}.$$
            Then, with the left action, the entries of $A^{-1}dA=A^*dA$ are invariant $1$-forms. We calculate
            
$$A^*dA=-\frac{dr}{r}I+\frac{1}{r^2}\pmatrix {\bar z_1 dz_1+\bar z_2 dz_2& -\bar z_1 d\bar z_2+\bar z_2 d\bar z_1\cr
           z_1 d z_2- z_2 d z_1 &  z_1 d\bar z_1+ z_2 d \bar z_2}=\pmatrix {i\sigma_1& -\sigma_2-i\sigma_3\cr
           \sigma_2+i\sigma_3 & -i\sigma_1}$$
  where
  \begin{eqnarray*}
  v_1&=&r^{-1}dr+i\sigma_1=(\bar z_1 dz_1+\bar z_2 dz_2)/r^2\\
  v_2&=&\sigma_2+i\sigma_3= (z_1 d z_2- z_2 d z_1)/r^2
  \end{eqnarray*}
  
  and these give a basis for the $(1,0)$-forms. We see that 
  $2dr=r(v_1+\bar v_1)$ so that
 $$\partial r=rv_1/2,\quad \bar\partial r=r\bar v_1/2$$ and hence
 $$\partial v_1=0,\quad \bar\partial v_1= dv_1=id\sigma_1=2i\sigma_2\sigma_3=-v_2\bar v_2.$$
Furthermore 
$$\partial v_2=v_1v_2,\quad \bar\partial v_2=-\bar v_1 v_2,
\quad \bar\partial(v_1v_2)=v_1\bar v_1v_2.$$
\vskip .25cm
We look for invariant solutions to the generalized K\"ahler equations where 
$$\rho_1=\exp \beta_1=\exp[dz_1dz_2]=\exp[r^2v_1v_2]$$
 and $\rho_2=\exp \beta_2$ where    
$$\beta_2=\sum_{i,j} H_{ij}v_i\bar v_j+ \lambda v_1v_2+\mu \bar v_1\bar v_2$$ 
(with $H_{ij},\lambda$ and $\mu$ functions of $r$) is a general invariant $2$-form.    
The algebraic compatibility conditions from Lemma 1 are:
$$(\beta_2-\beta_1)^2=0=(\beta_2-\bar\beta_1)^2$$   
which gives on subtraction
$$\beta_2(v_1v_2-\bar v_1\bar v_2)=0$$
or equivalently $\lambda=\mu$. We then get
$$0=\beta_2^2-2\beta_2\beta_1=(\sum H_{ij}v_i\bar v_j)^2+2\lambda^2v_1v_2\bar v_1\bar v_2-2\lambda r^2v_1v_2\bar v_1\bar v_2$$
or equivalently
\begin{equation}
\det H=\lambda(\lambda-r^2)
\label{detH}
\end{equation}
\vskip .5cm
We also know that $d\beta_2=0$ so that
\begin{eqnarray*}
\bar\partial (H_{ij}v_i\bar v_j)+\partial \lambda \bar v_1\bar v_2+\lambda \partial (\bar v_1\bar v_2)&=&0\\
\partial (H_{ij}v_i\bar v_j)+\bar\partial \lambda  v_1 v_2+\lambda \bar\partial  (v_1v_2)&=&0
\end{eqnarray*}
But $H$ and $\lambda$ are functions of $r$ and so from the first equation, expanding and collecting terms in $v_1\bar v_1\bar v_2$ we obtain
\begin{equation}
rH_{12}'+2H_{12}=r\lambda'-2\lambda
\label{H12}
\end{equation}
while collecting terms in $\bar v_1 v_2\bar v_2$ yields
\begin{equation}
rH_{22}'=2H_{11}
\label{H22}
\end{equation}
The second equation gives (\ref{H22}) again and also
\begin{equation}
rH_{21}'+2H_{21}=-r\lambda'+2\lambda
\label{H21}
\end{equation}
We can solve these by quadratures: from (\ref{H12}) we get
\begin{eqnarray*}
r^2H_{12}&=&\int^r (s^2\lambda' -2s\lambda)ds=r^2\lambda-4\int_a^r s\lambda ds\\
r^2H_{21}&=&-r^2\lambda+4\int_{a'}^r s\lambda ds.
\end{eqnarray*}
If we set
$$L(r)=\int_a^r s\lambda ds$$
then $\lambda=L'/r$ and 
$$r^2H_{12}=rL'-4L,\quad 
r^2H_{21}= -rL'+4L+b$$
and then $\det H=\lambda(\lambda-r^2)$ gives
$$H_{11}H_{22}=8\frac{LL'}{r^3}-16\frac{L^2}{r^4}-L'r+b\frac{L'}{r^3}-4b\frac{L}{r^4}.$$
Substituting $rH_{22}'=2H_{11}$ from (\ref{H22}) and integrating by parts leads to
\begin{equation}
H_{22}^2=16\frac{L^2}{r^4}+4b\frac{L}{r^4}-4L+c
\label{Hform}
\end{equation}
Thus an arbitrary complex function $L$ and three constants of integration $a,b,c$ give the general solution to the equations. Note for comparison that an $SU(2)$-invariant K\"ahler metric involves one \emph{real} function of $r$ -- the invariant K\"ahler potential.
\vskip .25cm
There is a lot of choice here but to produce an example let us take for simplicity $a=a'=0$ so that $H_{12}=-H_{21}$ and therefore $b=0$, and take $c=0$ so that
\begin{equation}
H_{22}^2=16\frac{L^2}{r^4}-4L
\label{examp1}
\end{equation} 
Let $L$ be real, then so is $\lambda$ and $H_{12}$. If $L$ negative then  $H_{22}^2$ is positive from  (\ref{examp1}). This means that  $H_{22}$ is real, and hence from (\ref{H22})  so is  $H_{11}$. Choose the positive square root for $H_{22}$.

Now $\beta_2=\sum H_{ij}v_i\bar v_j+ \lambda v_1v_2+\mu \bar v_1\bar v_2$ and $\lambda$ and $H_{ij}$ are real and $H_{12}=-H_{21}$ so
\begin{equation}
\beta_2-\bar\beta_2= 2(H_{11}v_1\bar v_1+H_{22}v_2\bar v_2)
\label{ibeta2}
\end{equation}
and for this to be symplectic  $H_{11}$ and $H_{22}$ must be non-zero. To get a  generalized K\"ahler metric we need from (\ref{posit}) to have $(\beta_1-\bar\beta_1)(X,\bar X)$ definite on the space of vectors $X$ with $i_X(\beta_1-\beta_2)=0.$ If $\nu_1,\nu_2,\bar\nu_1,\bar\nu_2$ is the dual basis to $v_1,v_2,\bar v_1,\bar v_2$ then $X$ must be  a linear combination of 
\begin{equation}
\lambda\nu_1-H_{12}\bar\nu_1+H_{11}\bar \nu_2,\quad \lambda\nu_2-H_{22}\bar\nu_1-H_{12}\bar\nu_2.
\label{CPhol}
\end{equation}
Since $\beta_1-\bar\beta_1=r^2(v_1v_2-\bar v_1\bar v_2)$ this gives $(\beta_1-\bar\beta_1)(X,\bar X)$ relative to this basis as the Hermitian form 
$$\pmatrix{2r^2\lambda H_{11}& \cr & 2r^2\lambda H_{22}}$$ 
so we also need $H_{11}$ to be positive.
\vskip .25cm
Notice now the point we have reached: $H_{11}$ and $H_{22}$ must be positive, which means that
\begin{equation}
H_{11}v_1\bar v_1+H_{22}v_2\bar v_2
\label{kform}
\end{equation}
is a positive definite Hermitian form. Moreover $rH_{22}'=2H_{11}$, and this implies that the form is K\"ahler. In fact if $\phi(r)$ satisfies $H_{22}=r\phi'/2$, this is the K\"ahler metric $i\partial\bar\partial \phi$, with $\phi$ as a K\"ahler potential.

Thus each $SU(2)$-invariant K\"ahler metric defines canonically, through the functions $L,\lambda$ and $H_{12}$ defined in terms of $H_{22}$, an $SU(2)$-invariant generalized K\"ahler metric on $\C^2\setminus \{0\}$.

\begin{prp} If the K\"ahler metric (\ref{kform}) extends to $\CP^2$, so does the generalized K\"ahler structure.
\end{prp}

\begin{prf} Since $\beta_1^{-1}$ is a global holomorphic Poisson structure on $\CP^2$, we know that the generalized complex structure $J_1$ extends to the whole of $\CP^2$, so 
we only need to check that $\beta_2$  also extends. We begin at $r=0$, the origin in $\C^2$. Clearly $r^2=z_1\bar z_1+z_2\bar z_2$ is smooth on $\C^2$. We shall use the fact that if $f(r)$ extends to a smooth function on a neighbourhood of the origin in $\C^2$ then $f(r)=f(0)+r^2f_1(r)$ where $f_1$ is also a smooth function.

If $g$ is the K\"ahler metric and $X=r\partial/\partial r$ the Euler vector field, then $g(X,X)=H_{11}$ is smooth on $\C^2$ and vanishes at the origin so 
 $H_{11}=r^2f_1$ for smooth $f_1>0$. The volume form of $g$ is $r^{-1}H_{11}H_{22}dr\sigma_1\sigma_2\sigma_3$ and comparing with the Euclidean volume $r^3dr\sigma_1\sigma_2\sigma_3$ we see that $H_{22}=r^2f_2$ for $f_2>0$ smooth.

Equation (\ref{examp1}) gives
$$L=\frac{r^4}{8}\left[1-\sqrt{1+(4H_{22}^2/r^4)}\right]=\frac{r^4}{8}\left[1-\sqrt{1+4f_2^2}\right]$$
and so $L=r^4 f_3$ for $f_3$ smooth.
By definition,
$\lambda=L'/r=4r^2f_3+r^3f_3'=r^2f_4$
since for any smooth $f(r)$
$$rf'=\sum x_i\frac{\partial f}{\partial x_i}$$
which is smooth. (In fact since this expression also vanishes at $0$ we have $rf'=r^2g$ for $g$ smooth.) 
 Since $r^2v_1v_2=dz_1dz_2$,  this shows that the term $\lambda (v_1v_2+ \bar v_1\bar v_2)$ is smooth.

 Now $r^2H_{12}=rL'-4L=r^5f_3'$ so $H_{12}=r^2(rf_3')=r^4f_5$ for smooth $f_5$,  which means that $H_{12}v_1\bar v_2$  and $H_{21}v_2\bar v_1$ are smooth since $r^2v_1=\bar z_1 dz_1+\bar z_2 dz_2, r^2v_2= z_1 d z_2- z_2 d z_1$. Hence the form $\beta_2$ is smooth at the origin.

 From (\ref{ibeta2}) the imaginary part of $\beta_2$ is nondegenerate at the origin since the K\"ahler metric is.
\vskip .25cm
As $r\rightarrow \infty$ we need to take homogeneous coordinates on $\CP^2$ so that $\C^2$ is parametrized by 
$[1,z_1,z_2]=[1/z_1,1,z_2/z_1]$, so we use  local affine coordinates  $w_1,w_2$ where for $z_1\ne 0$, 
$$w_1=\frac{1}{z_1},\quad w_2=\frac{z_2}{z_1}.$$
The projective line at infinity is then $w_1=0$. In these coordinates we have 
$$r^2=\frac{1+\vert w_2\vert^2}{\vert w_1\vert^2}$$
so $1/r^2$ is smooth and
\begin{equation}
v_1=\frac{\bar w_2dw_2}{1+\vert w_2\vert^2}-\frac{dw_1}{w_1},\quad v_2=\frac{\bar w_1dw_2}{w_1(1+\vert w_2\vert^2)}
\label{v1v2}
\end{equation}
Note here that
$$\frac{1}{r^2}v_1=\frac{\vert w_1\vert^2\bar w_2 dw_2}{(1+\vert w_2\vert^2)^2}-\frac{\bar w_1dw_1}{1+\vert w_2\vert^2}$$
is smooth  at $r=\infty$, and similarly $v_2/r^2, v_1\bar v_2/r^2$ are smooth. 

The coefficient of $dw_1d\bar w_1$ in $H_{11}v_1\bar v_1+H_{22}v_2\bar v_2$ is
$H_{11}/\vert w_1\vert^2$ so this is smooth and hence $r^2H_{11}=g_1$, a smooth function. Considering the coefficient of $dw_2d\bar w_2$ we see that $H_{22}$ is smooth.

 Now
\begin{equation}
L=\frac{r^4}{8}\left[1-\sqrt{1+4H_{22}^2/r^4}\right]=-\frac{1}{2}\frac{H_{22}^2}{1+\sqrt{1+4H_{22}^2/r^4}}
\label{ell}
\end{equation}
which is smooth and
$g_1=r^2H_{11}=r^3H_{22}'/2$ 
so that differentiating (\ref{ell}) 
$\lambda={L'}/{r}=g_2/r^4$ 
where $g_2$ is smooth. This means from (\ref{v1v2}) that $\lambda(v_1\bar v_1+v_2\bar v_2)$ is smooth. Finally
$H_{12}=\lambda-{4}L/{r^2}=g_3/r^2$
where $g_3$ is smooth, and so $H_{12}v_1\bar v_2$ is smooth. Thus $\beta_2$ extends as $r\rightarrow \infty$.

The argument  for $z_2\ne 0$ is the similar.
\end{prf}

\subsection{Example:  the Hirzebruch surface ${\mathbf F_2}$}\label{f2}

We can apply the above formalism with different boundary conditions to the Hirzebruch surface ${\mathbf F}_2$. Recall that this is 
$${\mathbf F}_2=P({\mathcal O}\oplus {\mathcal O}(-2))=P({\mathcal O}\oplus K)$$
since the canonical bundle $K$ of $\CP^1$ is ${\mathcal O}(-2)$.  The canonical symplectic form on $K$ extends to a meromorphic form $\beta_1$ on ${\mathbf F}_2$, and its inverse, a Poisson structure, defines the generalized complex structure $J_1$.
\vskip .25cm
On $K$ we take local coordinates $(w,z)\mapsto wdz$ where $z$ is an affine coordinate on $\CP^1$. Then for each quadratic polynomial $q(z)$ 
$$q(z)\frac{d}{dz}$$ 
is a global holomorphic vector field on $\CP^1$ so that
$$(w,z)\mapsto (w,wz,wz^2)$$
is a well defined map from $K$ to the cone $x_2^2=x_1x_3$ in $\C^3$. The map 
$$(z_1,z_2)\mapsto (z_1^2,z_1z_2,z_2^2)$$
maps the quotient $\C^2/\pm 1$ isomorphically to this cone and the  Hirzebruch surface is a compactification of the surface obtained by resolving  the singularity at the origin of this cone. Our ansatz above for $\C^2\setminus\{0\}$ extends to the quotient which is $\R^+\times SO(3)$ since we were using left-invariant forms. We need to adapt in a  different way  to extend at $r\rightarrow 0$ which is a rational curve of self-intersection $-2$ and $r\rightarrow \infty$, a rational curve of self-intersection $+2$.
\vskip .25cm
To proceed as $r\rightarrow 0$ we change coordinates from $z_1,z_2$ to $w,z$:
$$w=z_1^2,\qquad z=z_2/z_1.$$
Then $dzdw=2dz_2dz_1$, so here we see that the standard $2$-form on $\C^2$ is a multiple of the canonical symplectic form on the holomorphic cotangent bundle. We find
$$r^2=\vert w\vert(1+\vert z\vert^2)$$
so in particular $r^4$ is smooth. Furthermore
\begin{equation}
2\frac{dr}{r}=\frac{1}{2}\left[\frac{dw}{w}+\frac{d\bar w}{\bar w}\right]+\frac{d(z\bar z)}{1+\vert z\vert^2}
\label{diffr}
\end{equation}

We also calculate
\begin{equation}
v_1=\frac{dw}{2w}+\frac{\bar z dz}{1+\vert z\vert^2},\quad v_2=\frac{wdz}{\vert w\vert (1+\vert z\vert^2)}
\label{v12}
\end{equation}
Thus $r^2v_1v_2$ and $r^4v_1\bar v_1$ are smooth, and
$$v_2\bar v_2=\frac{dz d\bar z}{(1+\vert z\vert^2)^2}$$
which is smooth.
\vskip .25cm
Suppose in this case that $H_{11}v_1\bar v_1+H_{22}v_2\bar v_2$ extends as a K\"ahler form. Then considering the coefficient of $dwd\bar w$, $H_{11}=r^4 f_1$ where $f_1>0$ is smooth and $H_{22}$ itself is smooth and positive. The reality conditions on $H_{ij}$ are the same as the $\CP^2$ case and the constants of integration $a,a',b$ vanish as before but we now take $c$ in (\ref{Hform}) to be the limiting value $H_{22}^2(0)$. Since $H_{22}'>0$, $H_{22}^2-c>0$.
From (\ref{Hform}) we obtain
\begin{equation}
L=\frac{r^4}{8}\left[1-\sqrt{1+4(H_{22}^2-c)/r^4}\right].
\label{Lformula}
\end{equation}
We now use the familiar formula for a                                                                                                                      smooth function $f$ 
\begin{equation}
f(x)-f(x_0)=\sum_i(x-x_0)_i\int_0^1\frac{\partial f}{\partial x_i}(x_0+t(x-x_0))dt
\label{expand}
\end{equation}
where the coordinates $x_i$  are the real and imaginary parts of $z,w$ and we take $x_0=(z_0,0)$.   From  (\ref{diffr}) we calculate the derivatives 
$$\frac{\partial H_{22}}{\partial w}=\frac{r}{4w}H_{22}'=\frac{1}{2w}H_{11}=\frac{r^4 f_1}{2w}=\frac{1}{2}\bar w(1+\vert z\vert^2)^2f_1$$
(since $rH_{22}'=2H_{11}$ and $H_{11}=r^4f_1$ for smooth $f_1$) and 
$$\frac{\partial H_{22}}{\partial z}=\frac{r\bar z}{2(1+\vert z\vert^2)}H_{22}'=\frac{w\bar w f_1}{1+\vert z\vert^2}.$$
Putting these and their conjugates into (\ref{expand}) with $f=H_{22}$ we see that  $H_{22}(x)-H_{22}(x_0)=w\bar w f_2$
for a smooth function $f_2$ and hence from the formula for $L$ above 
$L=r^4f_3$ where $f_3$ is smooth. This gives 
$$\lambda=\frac{L'}{r}=4r^2f_3+4r^2w\frac{\partial f_3}{\partial w}.$$ This is $r^2f_4$ where $f_4$ is smooth and so $\lambda(v_1v_2+\bar v_1\bar v_2)$ is smooth since $r^2v_1 v_2$ is smooth.
Now 
\begin{equation}
H_{12}=\frac{L'}{r}-4\frac{L}{r^2}=r^3f_3'=4wr^2\frac{\partial f_3}{\partial w}=4w\vert w\vert (1+\vert z\vert^2)\frac{\partial f_3}{\partial w}
\label{h12}
\end{equation}
From (\ref{v12}) we see that $H_{12}v_1\bar v_2$ is smooth.
\vskip .25cm
In a neighbourhood of the curve $r=\infty$ we have coordinates $w'=1/w, z'=z$ and the calculations are very similar. In particular $1/r^4$ is smooth and $H_{22}$ is smooth and nonzero at infinity. Let $c'=\lim_{r\rightarrow\infty}H_{22}^2$. Then from (\ref{Lformula}) we have
$$L=-\frac{1}{4}(c'-c)+\frac{1}{r^4}g$$
where $g$ is smooth. This gives the required behaviour of $L$ and $\lambda$ for $\beta_2$ to extend to the curve at infinity.
\section{Bihermitian metrics}

\subsection{Generalized K\"ahler and bihermitian structures}

The generalized K\"ahler structures described above have a very concrete Riemannian description, owing to the following remarkable theorem of Gualtieri \cite{Gu}:

\begin{thm} \label{bi}
A generalized K\"ahler structure on a manifold $M^{2m}$ is equivalent to:
\begin{itemize}
\item
a Riemannian metric $g$
\item
two integrable complex structures $I_+,I_-$ compatible with the metric
\item
a $2$-form $b$ such that
$d^c_-\omega_-=db=-d^c_+\omega_+$ 
\end{itemize}
where $\omega_+,\omega_-$ are the two hermitian forms and $d^c=I^{-1}dI=i(\bar\partial-\partial)$.
\end{thm}

An equivalent description is to say that there are two connections $\nabla^+,\nabla^-$ which preserve the metric and the complex structures $I_+,I_-$ respectively and these are related to the Levi-Civita connection $\nabla$ by
\begin{equation}
\nabla^{\pm}=\nabla\pm\frac{1}{2}g^{-1}h
\label{del}
\end{equation}
where $h=db$ is of type $(2,1)+(1,2)$ with respect to both complex structures. In the K\"ahler case $I_+=I,I_-=-I$ and $b=0$. 

This is the geometry introduced $20$ years ago in the physics literature \cite{R} and more recently studied by differential geometers in four dimensions as ``bihermitian metrics", as in \cite{AGG}.
\vskip.25cm
Following \cite{Gu}, to derive this data from the generalized K\"ahler structure one looks at the  eigenspaces of $J_1J_2$. Since $J_1$ and $J_2$ commute, $(J_1J_2)^2=(-1)^2=1$. As before we choose $V_+$ to be the subbundle where $J_1=J_2$ and $V_-$ where $J_1=-J_2$. If the quadratic form  $(J_1J_2A,A)$ is negative definite,  the natural inner product on $T\oplus T^*$ is positive definite on  $V_+$, and negative definite on the complementary  eigenspace $V_-$. Since the signature of the quadratic form is $(2m,2m)$ each such space is $2m$-dimensional. Moreover since $T$ and $T^*$ are isotropic, $V_+\cap T=0=V_+\cap T^*$ and so $V_+$ is the graph of an invertible map from $T$ to $T^*$, i.e. a section $g+b$ of $T^*\otimes T^*$, where $g$ is the symmetric part and $b$ the skew-symmetric part. The bundle $V_+$ is preserved by $J_1$ and identified with $T$ by projection, and hence $J_1$ (or equivalently $J_2$) induces a complex structure $I_+$. Similarly on $V_-$, $J_1$ or $-J_2$ gives $T$ the complex structure $I_-$.

Conversely, as Gualtieri shows, given the bihermitian data above, the two commuting generalized complex structures are defined by 

\begin{equation}
J_{1/2}=\frac{1}{2}\pmatrix{1&0\cr
                 b&1}\pmatrix{I_+\pm I_- & -(\omega_+^{-1}\pm\omega_-^{-1})\cr
                              \omega_+\pm \omega_- & -(I_+^*\pm I_-^*)}\pmatrix{1&0\cr
                 -b&1}
\label{J12}
\end{equation}
\vskip .25cm
Our standard examples are constructed from closed forms $\rho_1=\exp\beta_1,\rho_2=\exp\beta_2$, so we look next at how the bihermitian structure is encoded in these.

The identification of $T$ with $V_+$ can be written as $X\mapsto X+(g(X,-)+b(X,-))$. If $X$ is a $(1,0)$-vector with respect to $I_+$ then this is 
$$X+\xi= X+i_X(b- i\omega_+)$$
where $\omega_+$ is the Hermitian form for $I_+$. If this lies in $E_1$, it annihilates $\exp \beta_1$, so $i_X\beta_1+ i_X(b- i\omega_+)=0$. Thus $\beta_1+b-i\omega_+$ is of type $(0,2)$, and similarly for  $E_2$. Thus there are $(0,2)$-forms $\gamma_1,\gamma_2$ such that  
$$\beta_1=-b+i\omega_+ +\gamma_1,\quad \beta_2=-b+i\omega_+ +\gamma_2.$$
Since $\beta_1,\beta_2$ are closed this means that $\gamma=\bar\beta_1-\bar\beta_2=\bar\gamma_1-\bar\gamma_2$ is a holomorphic $(2,0)$-form with respect to $I_+$.
The form $(\beta_1-\beta_2)$ defines the complex structure $I_+$ -- the $(1,0)$ vectors are the solutions to $i_X(\beta_1-\beta_2)=0$  and the metric on such $(1,0)$-vectors is given by $(\beta_1-\bar\beta_1)(X,\bar X)$. 

Changing to $V_-$,  the identification with $T$ is $X\mapsto X+(-g(X,-)+b(X,-))$ and then $\beta_1=-b+i\omega_- +\delta_1,\quad \bar\beta_2=-b+i\omega_- +\delta_2.$
where $\delta_1,\delta_2$ are $(0,2)$-forms with respect to $I_-$.
\vskip .25cm
 In four real dimensions we now give the precise relationship between the bihermitian description and the generalized K\"ahler one. First note that $\omega_-^{1,1}$ is self-dual and type $(1,1)$ so there is a real smooth function $p$ such that 
$$\omega_-^{1,1}=p\omega_+.$$
Moreover since $\omega_+^2=\omega_-^2$,  $\vert p\vert \le 1$.
 
From above we have
\begin{eqnarray*}
\beta_1&=&-b+i\omega_++\gamma_1=-b+i\omega_-+\delta_1\\
\beta_2&=&-b+i\omega_++\gamma_2=-b-i\omega_-+\bar\delta_2
\end{eqnarray*}
where $\gamma_1,\gamma_2$ are $(0,2)$ with respect to $I_+$ and $\delta_1,\delta_2$ are $(0,2)$ with respect to $I_-$. We let $\bar\gamma=\gamma_1-\gamma_2$ be the closed $(0,2)$ form, non-vanishing  since $\beta_1-\beta_2$ is non-zero from Lemma 1.

\begin{prp} \label{4formulas}In the terminology above,
\begin{itemize}
\item
$\beta_1=b+i\omega_+-(p-1)\bar\gamma/2$
\item
$\beta_2=b+i\omega_+-(p+1)\bar\gamma/2$
\item
$\omega_-=p\omega_++i(p^2-1)\bar\gamma/4-i(p^2-1)\gamma/4$
\end{itemize}
\end{prp}

\begin{prf}  Since we are in  two complex dimensions, there are functions $q_1,q_2$ such that the $(0,2)$ forms $\gamma_1,\gamma_2$ are given by  $\gamma_1=q_1\bar\gamma, \gamma_2=q_2\bar\gamma$ and since 
$$\beta_1-\beta_2=\gamma_1-\gamma_2=\bar\gamma$$
we have $q_1-q_2=1$. Similarly $\omega_-^{0,2}=r\bar\gamma$.

We have $\omega_-^2=\omega_+^2$ since this is the Riemannian volume form and 
$\omega_-=p\omega_++r\gamma+\bar r\bar\gamma$ since it is self-dual, hence
$$\omega_+^2=\omega_-^2=(p\omega_++r\gamma+\bar r\bar\gamma)^2=p^2\omega_+^2+2\vert r\vert^2\gamma\bar\gamma$$
and so
\begin{equation}
(1-p^2)\omega_+^2=2\vert r\vert^2\gamma\bar\gamma.
\label{eqA}
\end{equation}
Also $i\omega_++\gamma_1=i\omega_-+\delta_1$ and $\delta_1^2=0$ since it is of type $(0,2)$ relative to $I_-$ so 
$$0=(i\omega_++\gamma_1-i\omega_-)^2=(i\omega_++q_1\gamma-i[p\omega_++r\gamma+\bar r\bar\gamma])^2$$
and this gives
\begin{equation}
-(1-p)^2\omega_+^2=2i(q_1-ir)\bar r\gamma\bar\gamma.
\label{eqB}
\end{equation}
The same argument for $\delta_2$ gives
\begin{equation}
(1+p)^2\omega_+^2=2i(q_2+ir)\bar r\gamma\bar\gamma.
\label{eqC}
\end{equation}
From (\ref{eqA}),(\ref{eqB}),(\ref{eqC}) we obtain
$$q_1=\frac{2ir}{p+1},\quad q_2=\frac{2ir}{p-1}$$
and from $q_1-q_2=1$ it follows that $r=i(p^2-1)/4$ and hence 
$q_1=-(p-1)/2$ and $q_2=-(p+1)/2$.
\end{prf}
\begin{rmk} The function $p$ (which figures prominently as the \emph{angle function} in the calculations of \cite{AGG}) can be read off from the $2$-forms $\beta_1,\beta_2$ using the above formulas. Recall that the imaginary part of $\beta$ must be symplectic to define a generalized complex structure. We calculate the two Liouville volume forms:
$$(\beta_1-\bar\beta_1)^2=(p-1)\gamma\bar\gamma\qquad (\beta_2-\bar\beta_2)^2=-(p+1)\gamma\bar\gamma.$$
\end{rmk}

\subsection{Examples}

 Because of Theorem \ref{bi}, the constructions in (\ref{cp2}) and (\ref{f2}) using generalized complex structures  furnish us with bihermitian metrics. We now write these down. The complex structures $I_+,I_-$ are determined by the respective $(0,2)$ forms $\beta_1-\beta_2$ and $\beta_1-\bar\beta_2$. It is straightforward to see that
\begin{eqnarray*}
\lambda(\beta_1-\beta_2)&=& (H_{12}v_1+H_{22}v_2+\lambda\bar  v_1)(-H_{11}v_1+H_{12}v_2+\lambda \bar v_2)\\
\lambda(\beta_1-\bar\beta_2)&=& (H_{12}\bar v_1-H_{22}\bar v_2+\lambda v_1)(H_{11}\bar v_1+H_{12}\bar v_2+\lambda  v_2)
\end{eqnarray*}
The metric is obtained from the Hermitian form $\beta_1-\bar\beta_1$ on $(1,0)$ vectors. Using the basis of $(0,1)$ forms for $I_+$ given by the decomposition of $\beta_1-\beta_2$ above this turns out to be diagonal and the metric itself written as  
$$H_{11}\left[\frac{dr^2}{r^2-2\lambda + 2H_{12}}+\frac{r^2\sigma_1^2}{r^2-2\lambda-2H_{12}}\right]+H_{22}\left[\frac{r^2\sigma_2^2}{r^2-2\lambda + 2H_{12}}+\frac{r^2\sigma_3^2}{r^2-2\lambda-2H_{12}}\right].$$

\begin{rmk} If we replace the Poisson structure $\sigma=\partial/\partial z_1\wedge \partial/\partial z_2$ in our examples on $\CP^2$ or $ {\mathbf F}_2$ by $t\sigma$, then as $t\rightarrow 0$ the limiting generalized complex structure $J_1$ arises from a complex structure  and we should obtain simply a K\"ahler metric. 

This is equivalent to replacing the $2$-form $\beta_1$ by $t^{-1}\beta_1$. The differential equations for $H_{ij}$ remain the same but the algebraic constraint $\det H=\lambda(\lambda-r^2)$ becomes $\det H=\lambda(\lambda-t^{-1}r^2)$. The metric then becomes
$$tH_{11}\left[\frac{dr^2}{r^2-2t\lambda + 2tH_{12}}+\frac{r^2\sigma_1^2}{r^2-2t\lambda-2tH_{12}}\right]+tH_{22}\left[\frac{r^2\sigma_2^2}{r^2-2t\lambda + 2tH_{12}}+\frac{r^2\sigma_3^2}{r^2-2t\lambda-2tH_{12}}\right]$$
and removing the overall factor of $t$ this tends to the K\"ahler metric $H_{11}v_1\bar v_1+H_{22}v_2\bar v_2$ we started our constructions with.
\end{rmk}

 Concerning our examples of $\CP^2$ and ${\mathbf F}_2$, one should be careful to distinguish the various complex structures. In each case we took a complex structure which had a holomorphic Poisson structure and used that to define a \emph{generalized} complex structure $J_1$. We then found a generalized complex structure $J_2$ commuting with it and reinterpreted the pair as a bihermitian metric with two integrable complex structures $I_+$ and $I_-$.

It is well-known that $\CP^2$ has a unique complex structure so that all three complex structures are equivalent by a diffeomorphism in that case. However, all the Hirzebruch surfaces ${\mathbf F}_{2m}$ are diffeomorphic to $S^2\times S^2$. For $m> 0$ there is a unique holomorphic $SL(2,\C)$ action which has two   orbits of complex dimension one: a curve of self-intersection $+2m$ and one of $-2m$. 

The complex structures $I_+,I_-$ that arose from our construction admit a holomorphic $SU(2)$ action and  there are two spherical orbits of real dimension $2$ corresponding to $r=0$ and $r=\infty$.  We shall show that the sphere $S_0$ given by $r=0$ is not holomorphic with respect to $I_+$. 

Note first that the $2$-form $\beta_1=r^2v_1v_2$ vanishes on $S_0$ because $\beta_1$ has type $(2,0)$ in the ${\mathbf F}_2$ complex structure and $S_0$ is  holomorphic. Since $\lambda=r^2f_4$, this means that $\lambda(v_1v_2+\bar v_1 \bar v_2)$ vanishes on $S_0$. But from (\ref{h12}) 
$$H_{12}=4w\vert w\vert (1+\vert z\vert^2)\frac{\partial f_3}{\partial w}$$
so that
$$H_{12}v_1\bar v_2=4\frac{\partial f_3}{\partial w}(\bar w dw d\bar z +\bar z\vert w\vert ^2dzd\bar z)$$
and this vanishes on $S_0$ since $w=0$ there. Thus, restricted to $S_0$, all the terms in 
$\beta_1-\beta_2$ except $H_{11}v_1\bar v_1+H_{22}v_2 \bar v_2$ vanish,
and the latter is non-zero since it is the K\"ahler metric we started from.
However $\beta_1-\beta_2$ is a $(0,2)$-form in the complex structure $I_+$ and this must vanish on $S_0$ if it is a holomorphic curve.
\vskip .25cm
We conclude that, with the complex structure $I_+$ this must be the Hirzebruch surface ${\mathbf F}_0=\CP^1\times \CP^1$.

 \subsection{Holomorphic Poisson structures} \label {holp}
 
 Apostolov et al. in \cite{AGG} considered the four-dimensional bihermitian case where $I_+$ and $I_-$ define the same orientation and proved that the subset on which $I_+=\pm I_-$ is an anticanonical divisor with respect to both complex structures.  Now an anticanonical divisor is a holomorphic section of $\Lambda^2T^{1,0}$ -- a holomorphic bivector $\sigma$. Since $[\sigma,\sigma]$ is a holomorphic section of $\Lambda^3 T^{1,0}$, in two complex dimensions this 
automatically vanishes and we have a Poisson structure. All compact surfaces with holomorphic Poisson structure have been listed by Bartocci and Macr\`\i  \,\,using the classification of complex surfaces \cite{Bart}, so considering this list  provides a basis for seeking compact bihermitian metrics in this dimension. Particular cases (overlooked in \cite{AGG}) are the projective bundle $P(1\oplus K)$ over any compact algebraic curve $C$, and the ``twisted' version $P(V)$ where 
$$0\rightarrow K\rightarrow V\rightarrow 1\rightarrow 0$$
is the nontrivial extension in $H^1(C,K)\cong \C$. When $C=\CP^1$ these two surfaces are are ${\mathbf F}_2$ and ${\mathbf F}_0$ respectively.

We show  now that the Poisson structure appears naturally in higher dimensions too.
\vskip .25cm
Let $M$ be a generalized K\"ahler manifold, now considered from the bihermitian point of view. Following \cite{AGG} we consider the $2$-form
$$S(X,Y)=g([I_+,I_-]X,Y).$$
Since
\begin{eqnarray*}
S(I_+X,I_+Y)&=&g(I_+I_-I_+X,I_+Y)-g(I_-I_+^2X,I_+Y)\\
&=&g(I_-I_+X,Y)+g(I_-X,I_+Y)\\
&=&=g([I_-,I_+]X,Y)=-S(X,Y)
\end{eqnarray*}
this form is of type $(2,0)+(0,2)$.
Pick the complex structure $I_+$. Using the antilinear isomorphism $T^{1,0}\cong (\bar T^*)^{0,1}$ provided by the hermitian metric, its $(0,2)$ part can be identified  with a section $\sigma_+$ of the bundle $\Lambda^2T^{1,0}$.

\begin{prp} \label{biv} The bivector $\sigma_+$ is a holomorphic Poisson structure.
\end{prp}

\begin{prf} We shall first show that $\sigma_+$ is holomorphic, and then that its Schouten bracket vanishes. 

Let $z_1,\dots,z_n$ be local holomorphic coordinates, then 
$$\sigma_+=\sum (I_- dz_i,dz_j)\frac{\partial}{\partial  z_i}\wedge\frac{\partial}{\partial  z_j}$$
where we use the inner product on $1$-forms defined by the metric and the complex structure $I_-$ on $1$-forms. We need to show that the functions $(I_- dz_i,dz_j)$ are holomorphic. Now
\begin{equation}
\frac{\partial}{\partial \bar z_k}(I_- dz_i,dz_j)=((\nabla^+_{\bar k}I_-) dz_i,dz_j)+(I_-\nabla^+_{\bar k}dz_i,dz_j)+(I_- dz_i, \nabla^+_{\bar k}dz_j).
\label{holo}
\end{equation}
The Levi-Civita connection $\nabla$ has zero torsion so 
$$0=d(dz_i)=\sum dz_k\wedge \nabla_k dz_i+\sum d\bar z_k \wedge \nabla_{\bar k}dz_i.$$
But from (\ref{del}) $\nabla^+=\nabla +H/2$ where $H=g^{-1}db$, so
\begin{equation}
0=\sum_k dz_k\wedge (\nabla^+_k-H_{k}/2) dz_i+\sum_k d\bar z_k \wedge (\nabla^+_{\bar k}-H_{\bar k}/2)dz_i.
\label{dbareq}
\end{equation}
Now $\nabla^+$ preserves $I_+$ so that $\nabla^+_k dz_i$ and $\nabla^+_{\bar k} dz_i$ are $(1,0)$-forms. However, since $H$ is of type $(2,1)+(1,2)$, $H_k(dz_i)$ has a $(0,1)$ component. Equating the $(1,1)$ component of (\ref{dbareq}) to zero, the two contributions of $H$ give 
\begin{equation}
\nabla^+_{\bar k}dz_i=H_{\bar k}(dz_i)
\label{nab1}
\end{equation}
Now $I_-$ is preserved by $\nabla^-$ and from (\ref{del}) $\nabla^-=\nabla^+-H$, so
$$\nabla_{\bar k}^+ I_-=[H_{\bar k},I_-].$$
Using this and (\ref{nab1}) in (\ref{holo}) we obtain 
$$\frac{\partial}{\partial \bar z_k}(I_- dz_i,dz_j)=([H_{\bar k},I_-] dz_i,dz_j)+(I_-H_{\bar k}(dz_i),dz_j)+(I_- dz_i, H_{\bar k}(dz_j))=0$$
and so $\sigma_+$ is holomorphic.
\vskip .25cm
To prove that $\sigma_+$ is Poisson we use (\ref{J12}) and the observation that the upper triangular part of $J_1$  is a real Poisson structure. This means that
$$[\omega_+^{-1}+\omega_-^{-1},\omega_+^{-1}+\omega_-^{-1}]=0.$$
Now since $\omega_+$ is of type $(1,1)$, $\omega_+^{-1}+\omega_-^{-1}=h+\sigma_++\bar\sigma_+$ where $h$ is a bivector of type $(1,1)$. Because $\sigma_+$ is holomorphic, $[h,\sigma_+]$ has no $(3,0)$ component and so the $(3,0)$ component of $0=[h+\sigma_++\bar\sigma_+,h+\sigma_++\bar\sigma_+]$ is just $[\sigma_+,\sigma_+]$. Hence  $[\sigma_+,\sigma_+]=0$ and we have a holomorphic  Poisson structure.
\end{prf}
\vskip .25cm
When the generalized K\"ahler structure is defined by $\rho_1=\exp \beta_1,\rho_2=\exp \beta_2$, as in Lemma 1, $\sigma_+$ has a direct interpretation. Recall that $\bar\beta_1-\bar\beta_2 =\gamma$  is a
non-degenerate holomorphic $2$-form with respect to $I_+$. Then

\begin{prp} \label{gprop} Let $\sigma_+:(T^{1,0})^*\rightarrow T^{1,0}$ be the holomorphic Poisson structure corresponding to the  generalized K\"ahler structure given by $2$-forms $\beta_1,\beta_2$, and let $\gamma=\bar\beta_1-\bar\beta_2:T^{1,0}\rightarrow (T^{1,0})^*$ be the holomorphic $2$-form. Then 
$$\sigma_+=2i\gamma^{-1}.$$
\end{prp}

\begin{prf} From (\ref{J12}) $\sigma_+$ is given by the upper-triangular part of $J_1$ evaluated on one-forms of type $(1,0)$ with respect to $I_+$. Since $\gamma$ is a non-degenerate $(2,0)$ form, any $(1,0)$ form can be written  $i_X\gamma$ for a $(1,0)$-vector $X$. So we require to prove that if $X$ is a $(1,0)$ vector, then the $(1,0)$ component of 
$J_1(i_X\gamma)$ is $2iX$. Now 
\begin{eqnarray*}
J_1(i_X\gamma)&=&J_1(i_X(\bar\beta_1-\bar\beta_2))\\
&=& J_1(i_X\bar\beta_1-X+X-i_X\bar\beta_2)
\end{eqnarray*}
and by the definition of $J_1$, 
\begin{equation}
J_1(i_X\bar\beta_1-X)=-i(i_X\bar\beta_1-X)
\label{first}
\end{equation}
The term $X-i_X\bar\beta_2$ is acted on as $-i$ by $J_2$ and we split it into components for the two $J_1$ eigenspaces:
$$X-i_X\bar\beta_2=Y-i_Y\bar\beta_2+Z-i_Z\bar\beta_2.$$
Since $Z-i_Z\bar\beta_2$ is in the $-i$-eigenspace of both $J_1$ and $J_2$, $Z$ is of type $(0,1)$. Since $X=Y+Z$, $X=Y^{1,0}$. Now
$$J_1(X-i_X\bar\beta_2)=i(Y-i_Y\bar\beta_2)-i(Z-i_Z\bar\beta_2)$$
and adding this to (\ref{first}), the upper triangular part of $J_1$ is given by 
$$J_1(i_X(\bar\beta_1-\bar\beta_2))=2iX -2iZ$$
whose $(1,0)$ part is $2iX$.
\end{prf}
\begin{ex}
The examples of $\CP^2$ and ${\mathbf F}_2$ were constructed by using $2$-forms $\beta_1,\beta_2$. Since $\beta_1$ had a pole on the curve at $r=\infty$ and $\beta_2$ was smooth everywhere, the Poisson structures $\sigma_+=2i(\bar\beta_1-\bar\beta_2)^{-1}$ and $\sigma_-=2i(\bar\beta_1-\beta_2)^{-1}$ vanish there. 
\end{ex}

\section{Moduli spaces of instantons}

\subsection{Stability}
On a  $4$-manifold with a Hermitian structure, the anti-self-dual (ASD) $2$-forms are the $(1,1)$-forms orthogonal to the Hermitian form. Thus on a generalized K\"ahler $4$-manifold, a connection with anti-self-dual curvature (an instanton) has curvature of type $(1,1)$ with respect to both complex structures $I_+,I_-$. In fact, where $I_+\ne \pm I_-$, anti-self-duality is equivalent to this condition.

The equations $d^c_-\omega_-=db=-d^c_+\omega_+$ imply that $$dd_{\pm}^c\omega_{\pm}=0$$
which means that the metric is a Gauduchon metric with respect to both complex structures. With a Gauduchon metric one defines the \emph{degree} of a holomorphic line bundle $L$ by
$$\deg L=\frac{1}{2\pi}\int_M F\wedge\omega$$
where $F$ is the curvature form of a connection on $L$ defined by a Hermitian metric. Since a different choice of metric changes $F$ by $dd^cf$, the condition $dd^c\omega=0$ and integration by parts shows that the degree, a real number, is  independent of the choice of Hermitian metric on $L$. It has the usual property of degree that if a holomorphic section of $L$ vanishes on a divisor $D$ then 
$$\deg L=\int_D\omega.$$
So line bundles with sections which vanish somewhere have positive degree. 

\begin{rmk} Let us consider this non-K\"ahler degree for a  bihermitian surface  such that the Poisson structure vanishes on a divisor, like our examples of $\CP^2$ and $\CP^1\times \CP^1$, and assume for convenience that the surface also carries a K\"ahler metric. The canonical bundle $K$ has no holomorphic sections since the product with the Poisson structure, a section of $K^*$, would give a holomorphic function with zeroes. This means $H^{2,0}(M)=0$ and so $H^2(M)$ is purely of type $(1,1)$.

Now suppose that one of the generalized complex structures is defined by $\exp \beta$ where $\beta$ is closed. We saw in (\ref{4formulas}) that
 $\beta=-b+i\omega_++\gamma_1$ where $\gamma_1$ is of type $(0,2)$, so that the $(1,1)$ component of $\beta-\bar\beta$ is $2i\omega_+$. Thus the integral of $\omega_+$ over a holomorphic curve $C$, which is positive, is the same as the integral of the \emph{closed} form $(\beta-\bar\beta)/2i$. Let $W$ be the cohomology class of this form. Then we see that for every effective divisor $D$ on $M$, $WD>0$. Furthermore, $W$ is represented by the form 
$$(\beta-\bar\beta)/2i=\omega_+-i(\gamma_1-\bar\gamma_1)/2$$
which is self-dual, hence $W^2>0$. It follows from Nakai's criterion that $W$ is the cohomology class of a K\"ahler metric.

Since the ample cone generates the whole of the cohomology, we see that the non-K\"ahler degree in this case agrees with the ordinary K\"ahler degree of some K\"ahler metric. 

Observe also that $\beta-\bar\beta$ is also equal to $2i\omega_-+\delta_1-\bar\delta_1$ so that we obtain the same degree function on cohomology  for $I_+$ and $I_-$.
\end{rmk}
\vskip .25cm
Using this definition of degree, one can define the slope of a subbundle, and from that the stability of a holomorphic bundle. The key theorem in the area, proved by Buchdahl \cite{Buch} for surfaces and Li and Yau \cite{LY} in the general case, is that a bundle is stable if and only if it has an irreducible ASD connection. A good reference for this is the book \cite{LT}.

From this we already see that the moduli space ${\mathcal M}$ of ASD connections on a generalized K\"ahler manifold has two complex structures, by virtue of being the moduli space of stable bundles for both $I_+$ and $I_-$.

We shall prove the following theorem:

\begin{thm} \label{GKmod} Let $M^4$ be a compact even generalized K\"ahler manifold. Then the smooth points of the moduli space of ASD connections on a  principal $SU(k)$-bundle over $M$ carries a natural bihermitian metric such that $d^c_-\omega_-=H=-d^c_+\omega_+$ for some exact $3$-form $H$ of type $(2,1)+(1,2)$.
\end{thm}

From Gualtieri's theorem this has a generalized K\"ahler interpretation once we  choose a $2$-form $b$ such that $db=H$. 

\begin{rmk} In general, the moduli space of stable bundles may have singularities if  the obstruction space $H^2(M,\End_0 E)$ (where $\End_0$ denotes trace-free endomorphisms) is non-vanishing. However, if the Poisson structure $s$ on $M$ is non-zero, then
$$s:H^0(M,\End_0 E\otimes K)\rightarrow H^0(M,\End_0 E)$$
is injective. But stable bundles are simple,so  $H^0(M,\End_0 E)=0$. We deduce that  $H^0(M,\End_0 E\otimes K)$, and hence also its Serre dual $H^2(M,\End_0 E)$, must vanish, so the moduli space is smooth (see \cite{Bot1}). 

This vanishing also gives us by Riemann-Roch the dimension of the $SU(k)$ moduli space 
$$\dim_{\C}{\mathcal M}=2kc_2(E)-(k^2-1)\frac{1}{12}(c_1^2+c_2)(M).$$
The simplest case would be $k=2, c_2(E)=n$ for our examples $\CP^2,{\mathbf F}_2$ (or any rational surface) where $\dim_{\C}{\mathcal M}=4n-3.$
\end{rmk}

\subsection{The metric on the  moduli space}

In \cite{LT}  the metric structure of the moduli space of instantons on a Gauduchon manifold is discussed.  It differs in general from the Riemannian or K\"ahler case. In the Riemannian situation, the space of all connections is viewed as an infinite-dimensional affine space with group of translations $\Omega^1(M,\lie{g})$ and ${\mathcal L}^2$ metric
$$(a_1,a_2)=-\int_M\tr(a_1\wedge \star a_2).$$
 The solutions to the ASD equations form an infinite-dimensional submanifold with induced metric, and its quotient by the group of gauge transformations ${\mathcal G}$ is the moduli space, which acquires the quotient metric. To define this, one identifies the tangent space of the quotient at a point $[A]$ with the orthogonal complement to the tangent space of the gauge orbit at the connection $A$, with its  restricted inner product. The orthogonal complement is identified with the bundle-valued $1$-forms $a\in \Omega^1(M,\lie{g})$ which satisfy the equation
\begin{equation}
d_A^*a\, (=-\star  \dstar a)=0
\label{horizont}
\end{equation}
As the authors of \cite{LT} point out, this metric in the Gauduchon case is not Hermitian with respect to the natural complex structure that the moduli space acquires through its identification with the moduli space of stable bundles. Instead of the orthogonality (\ref{horizont}), one takes a different horizontal subspace defined by
\begin{equation}
\omega \wedge d_A^ca=0.
\label{newhorizont}
\end{equation}
\begin{lem} \label{hori}$\omega \wedge d^c_A a = \dstar a- d^c\omega \wedge a.$
\end{lem}
From this we see that when the metric is K\"ahler, $d^c\omega=0$, and so the two horizontality conditions coincide.
\begin{lemprf}
Note that for any $\psi\in \Omega^0(M,\lie{g})$,
\begin{equation}
d_A^c(\omega \wedge\tr(a \psi))=d^c\omega\wedge\tr(a\psi)+\omega\wedge\tr(d^c_A a\psi)-\omega\wedge\tr(a \wedge d_A^c\psi)
\label{exp1}
\end{equation}
and  $d_A^c\psi=I^{-1}d_AI\psi=-Id_A\psi$, so that $$\omega\wedge\tr(a\wedge d_A^c\psi)=-\omega\wedge\tr(a\wedge Id_A\psi)=(a,d_A\psi)\omega^2=\tr(\star a \wedge d_A\psi).$$
Integrating (\ref{exp1}) and using Stokes' theorem and the relation above, we get
$$\int_M[d^c\omega\wedge\tr(a\psi) +\omega\wedge\tr(d^c_A a\psi)-\tr(\dstar a \psi)]=0$$ so that
\begin{equation}
\omega \wedge d^c_A a = \dstar a- d^c\omega \wedge a.
\label{newhorizont1}
\end{equation}
\end{lemprf}

 With this choice of horizontal, the metric on the moduli space is Hermitian with Hermitian form
$$\tilde\omega(a_1,a_2)=\int_M\omega\wedge \tr(a_1\wedge a_2).$$
It is  shown in \cite{LT} that $\tilde\omega$ satisfies $dd^c\tilde\omega=0$.
The horizontal subspace (\ref{newhorizont}) defines a connection on the infinite-dimensional principal ${\mathcal G}$-bundle over the moduli space and  its curvature turns out to be of type $(1,1)$ on ${\mathcal M}$ (see \cite{LT}). We shall make use of these facts later.
\vskip .25cm
In order to prove Theorem \ref{GKmod} we need first to show that the application of L\"ubke and Teleman's approach to the two complex structures $I_+$ and $I_-$ yields the same metric.
\vskip .25cm
 The tangent space to the moduli space at a smooth point is the first cohomology of the complex:
$$\Omega^0(M,\lie{g})\stackrel{d_A}\longrightarrow \Omega^1(M,\lie{g})\stackrel{d^+_A}\longrightarrow \Omega_+^2(M,\lie{g})$$
where here the $+$ refers to projection onto the self-dual part. The metric is the induced inner product on the subspace of $\Omega^1(M,\lie{g})$ defined by the horizontality condition $\omega \wedge d^c_A a=0$. We shall write $[a]$ for the tangent vector to the moduli space represented by $a$.

In our case we have two such horizontality conditions $\omega_- \wedge d^c_-a=0$ and $\omega_+ \wedge d^c_+a=0$ (suppressing the subscript $A$ for clarity) and two representatives $a$ and $a+d_A\psi$ for the same tangent vector. We shall call these plus- and minus- horizontal respectively. We prove:
 
\begin{lem} Let $a$ and $a+d_A\psi$ satisfy 
$$\omega_-\wedge  d^c_-a=0,\quad \omega_+ \wedge d^c_+(a+d_A\psi)=0.$$
Then $(a,a)=(a+d_A\psi,a+d_A\psi)$.
\end{lem}

\begin{lemprf}
Since in our case $d^c_-\omega_-=db=h=-d^c_+\omega_+$ our two horizontality conditions are, from (\ref{newhorizont1}) 
$$ \dstar a- h \wedge a=0\quad \dstar\, (a+d_A\psi)+h \wedge (a+d_A\psi)=0$$ and so, eliminating $h\wedge a$,
$$2\dstar a+d_A\star d_A\psi+h\wedge d_A\psi=0.$$
This gives on integration
$$\int_M[2\tr(\dstar a \psi)+\tr(\dstar d_A\psi \psi)+h\wedge\tr(d_A\psi \psi)]=0.$$
But $\tr(d_A\psi\psi)=d\tr\psi^2/2$ so the last term is $d[(\tr\psi^2)h/2]$ as $h$ is closed. By Stokes' theorem we get 
$$2(a,d_A\psi)+(d_A\psi,d_A\psi)=0$$
and hence
$$(a+d_A\psi,a+d_A\psi)=(a,a)$$
as required.
\end{lemprf}
\subsection{The bihermitian structure}
So far, we have seen that 
${\mathcal M}$ has two complex structures
and a metric, Hermitian with respect to both.
We now need to show that 
$d^c_+\tilde\omega_+=H=-d^c_-\tilde \omega_-$ for an exact $3$-form $H$.
\vskip .25cm
Denote by ${\mathcal A}$ the affine space of all  connections on the principal bundle, then  a 
tangent vector is given by $a\in \Omega^1(M,\lie{g})$
and for any $2$-form $\omega$,
$$\Omega(a_1,a_2)=\int_M\omega\wedge \tr(a_1\wedge a_2)$$
is a closed and gauge-invariant $2$-form on ${\mathcal A}$. It is closed because it is translation-invariant on $\mathcal{A}$ (has ``constant coefficients").

We defined Hermitian forms $\tilde\omega_{\pm}$ on ${\mathcal M}$  by 
$$\tilde\omega_{\pm}([a_1],[a_2])=\Omega_{\pm}(a_1,a_2)=\int_M\omega_{\pm}\wedge\tr(a_1\wedge a_2)$$
where $a_1,a_2$ are plus/minus-{\it horizontal}. Now the formula for the exterior derivative of a $2$-form $\alpha$ is 
$$3d\alpha(a_1,a_2,a_3)=a_1\cdot\alpha(a_2,a_3)-\alpha([a_1,a_2],a_3)+\cyc$$
so, since $\Omega$ is closed 
$$3d\tilde\omega([a_1],[a_2],[a_3])=-\int_M\omega\wedge\tr([a_1,a_2]_V)\wedge a_3)+\cyc$$
where $[a_1,a_2]_V$ is the vertical component of the Lie bracket of the two vector fields. By definition this is the curvature of the ${\mathcal G}$-connection. If $\theta(a_1,a_2)\in \Omega^0(M,\lie{g})$ is this curvature then $[a_1,a_2]_V=d_A\theta(a_1,a_2)$. Using Stokes' theorem 
\begin{eqnarray*}
3d\tilde\omega([a_1],[a_2],[a_3])&=&-\int_M\omega\wedge\tr(d_A\theta(a_1,a_2)\wedge a_3)+\cyc\\
&=& \int_M d\omega\wedge \tr(\theta(a_1,a_2)a_3)+ \int_M\omega\wedge\tr(\theta(a_1,a_2)d_Aa_3)+\cyc\\
&=&\int_M d\omega\wedge \tr(\theta(a_1,a_2)a_3) +\cyc
\end{eqnarray*} 
since   $d_Aa_3$ is anti-self-dual and $\omega$ is self-dual so $\omega\wedge  d_Aa_3=0$.
\vskip .25cm
Now 
$d^c\omega(a_1,a_2,a_3)=-d\omega(Ia_1,Ia_2,Ia_3)$
and from \cite{LT} the 
curvature of the ${\mathcal G}$-bundle is of type $(1,1)$. This means that $\theta(Ia_2,Ia_3)=\theta(a_2,a_3)$ 
and so, for the structure $I_-$
\begin{equation}
d_-^c\tilde\omega([a_1],[a_2],[a_3])=\int_M d_-^c\omega_-\wedge  \tr(\theta_-(a_1,a_2)a_3)+\cyc
\label{deec}
\end{equation}
with a similar equation for $I_+$.
\vskip .25cm
To proceed further we need more information about the curvature $\theta(a_1,a_2)$.  On the affine space ${\mathcal A}$ the Lie bracket of $a_1$ and $a_2$ considered as vector fields is just 
$a_1\cdot a_2-a_2\cdot a_1$ where $a\cdot b$  denotes the flat derivative of $b$ in the direction $a$. The horizontality condition imposes a  constraint:
$$ \dstar a_2- h \wedge a_2=0.$$
Differentiating the constraint in the direction $a_1$ gives
$$[a_1,\star a_2]+\dstar a_1\cdot a_2-h\wedge a_1\cdot a_2=0.$$
The vertical component of the Lie bracket is $d_A\theta(a_1,a_2)$ which thus satisfies
\begin{equation}
\dstar d_A\theta-h \wedge d_A\theta +2[a_1,\star a_2]=0.
\label{thetaeq}
\end{equation}
Define the second order operator $\Delta:\Omega^0(M,\lie{g})\rightarrow \Omega^4(M,\lie{g})$ by
$$\Delta\psi=d_A\star d_A\psi-h\wedge d_A\psi,$$ then its formal adjoint is 
$$\Delta^*\psi=d_A\star d_A\psi+h\wedge d_A\psi$$
and  we rewrite (\ref{thetaeq}) as
\begin{equation}
\Delta\theta(a_1,a_2)+2[a_1,\star a_2]=0
\label{thetaeqs}
\end{equation}
for plus-horizontal vector fields $a_i$. Let $b_i=a_i+d_A\psi_i$ be the minus-horizontal representatives of $[a_i]$. By minus-horizontality we have 
$$0=d_A\star b_i+h\wedge b_i=\dstar\, (a_i+d_A\psi_i)+h\wedge  (a_i+d_A\psi_i)=\Delta^*\psi_i+\dstar a_i +h\wedge a_i$$
and together with the plus-horizontality condition  $ \dstar a_i- h \wedge a_i=0$ we get
\begin{equation}
2h\wedge a_i=-\Delta^*\psi_i.
\label{hai}
\end{equation}
Since $d^c_-\omega_-=h$, each integrand on the right hand side of (\ref{deec}) is, from  (\ref{hai}), of the form
$$h\wedge \tr(\theta(a_1,a_2)a_3)=-\tr(\theta(a_1,a_2)\Delta^*\psi_3/2).$$
 Performing the integration and using Stokes' theorem, we obtain
$$-\int_M \tr(\theta(a_1,a_2)\Delta^*\psi_3)/2=-\int_M \tr(\Delta\theta(a_1,a_2)\psi_3)/2=\int_M \tr([a_1,\star a_2]\psi_3)$$
from (\ref{thetaeqs}). 
\vskip.25cm
Working with the curvature of the plus-connection we get a similar expression so that we have two formulae:
\begin{eqnarray*}
d_-^c\tilde\omega_-([a_1],[a_2],[a_3])&=&\int_M \tr([a_1,\star a_2]\psi_3)+\cyc\\
d_+^c\tilde\omega_+([a_1],[a_2],[a_3])&=&-\int_M \tr([b_1,\star b_2]\psi_3)+\cyc
\end{eqnarray*}
Thus to 
obtain  $d_-^c\tilde\omega_- =-d_+^c\tilde\omega_+$, using $b_i=a_i+d_A\psi_i$  in the above leads to the need to prove:
\begin{lem}
$$\int_M [\tr([a_1,\star d_A\psi_2]\psi_3)+\tr([d_A\psi_1,\star a_2]\psi_3)+\tr([d_A\psi_1,\star d_A\psi_2]\psi_3)]+\cyc=0.$$
\end{lem}
\begin{lemprf}
Picking out the integrand involving $a_1$ in the cyclic sum we have 
\begin{eqnarray*}
\tr([a_1,\star d_A\psi_2]\psi_3)+\tr([d_A\psi_3,\star a_1]\psi_2)&=&\tr(\star a_1\wedge([\psi_3, d_A\psi_2]+ [d_A\psi_3,\psi_2]))\\
&=&-\tr(\star a_1\wedge d_A[\psi_2,\psi_3])
\end{eqnarray*}
and on integrating,
this is
\begin{eqnarray*}
-\int_M\tr(\star a_1\wedge d_A[\psi_2,\psi_3])&=&-\int_M\tr(\dstar a_1[\psi_2,\psi_3])\\
&=& -\int_M h\wedge \tr(a_1[\psi_2,\psi_3])\\
&=& \int_M\tr(\Delta^*\psi_1[\psi_2,\psi_3])/2
\end{eqnarray*}
from (\ref{hai}). But from the definition of $\Delta^*$ this is
$$\frac{1}{2}\int_M\tr(d_A\star d_A\psi_1[\psi_2,\psi_3])-\frac{1}{2}\int_M h\wedge \tr(d_A\psi_1[\psi_2,\psi_3]).$$
The cyclic sum of the second term vanishes  since
$$d\tr(\psi_1[\psi_2,\psi_3])=\tr(d_A\psi_1[\psi_2,\psi_3])+\cyc$$
and $h$ is closed.
Using Stokes' theorem on the first and expanding,  the cyclic sum   gives  
$$\frac{1}{2}\int_M\tr(\star d_A\psi_1\wedge( [d_A\psi_2,\psi_3]+[\psi_2,d_A\psi_3])+\cyc$$
which is 
$$-\int_M\tr([d_A\psi_1,\star d_A\psi_2]\psi_3)]+\cyc$$
and this proves the lemma.
\end{lemprf}
\vskip .25cm
We finally need to show that $H$ is exact. One might expect that we simply define a $2$-form $\tilde b$ from the $2$-form $b$ on $M$ by
\begin{equation}
\tilde b([a_1],[a_2])=\int_M b\wedge \tr(a_1\wedge a_2)
\label{btilde}
\end{equation}
to get $d\tilde b=d_-^c\tilde\omega_-$ but this does not hold.  The equation for the exterior derivative of $\tilde b$ gives  
$$3d\tilde b([a_1],[a_2],[a_3])=\int_M db\wedge \tr(\theta(a_1,a_2)a_3)+ \int_M b\wedge\tr(\theta(a_1,a_2)d_Aa_3)+\cyc.$$
When we used this above with $\omega_+, \omega_-$ replacing $b$, the second term vanished because $d_Aa_3$ is anti-self-dual and $\omega_{\pm}$ are self-dual. This is not the case for a general $b$, and will only be true if $b$ is self-dual. We shall see in Section 5 a more general occurrence of this phenomenon. However we do have the following:

\begin{lem} Any $2$-form $b$ on a compact oriented four-manifold  $M$ is the sum of a closed form and a self-dual form.
\end{lem}

\begin{lemprf} Use the non-degenerate pairing on $2$-forms
$$(\alpha,\beta)=\int_M\alpha\wedge \beta.$$ 
The annihilator of the self-dual forms $\Omega^2_+$ in this pairing is $\Omega^2_-$, and the annihilator  of $\Omega^2_{closed}$ is $\Omega^2_{exact}$  so the annihilator of $\Omega^2_++\Omega^2_{closed}$ is the intersection of $\Omega^2_-$ and $\Omega^2_{exact}$. But if $\alpha$ is exact, then by Stokes' theorem 
$$\int_M\alpha\wedge \alpha=0$$
and if $\alpha\in \Omega^2_-$ 
$$\int_M\alpha\wedge \alpha=-(\alpha,\alpha)$$
so if both hold then $\alpha=0$.
\end{lemprf}

It follows from this that $db=db_+$ where $b_+$ is self-dual, and then (\ref{btilde}) does define a form $\tilde b_+$ on the moduli space. It follows than that  $d\tilde b_+=d_-^c\tilde\omega_-=-d_+^c\tilde\omega_+$.

\subsection{The Poisson structures on ${\mathcal M}$}

As we saw in Proposition \ref{biv}, a generalized K\"ahler structure defines a holomorphic Poisson structure for each of the complex structures $I_+,I_-$. We shall determine these on  the instanton moduli  space next.

On the moduli space of stable bundles over a Poisson surface $M$, there is a canonical holomorphic Poisson structure, defined by Bottacin in \cite{Bot1}  as follows. The holomorphic tangent space at a bundle $E$ is the sheaf cohomology group $H^1(M,\End E)$ and by Serre duality, the cotangent space is $H^1(M,\End E \otimes K)$. The Poisson structure on $M$ is a holomorphic section $s$ of the anticanonical bundle $K^*$ and for $\alpha,\beta\in  H^1(M,\End E \otimes K)$, the Poisson structure $\sigma$ on the moduli space is defined by taking $\tr(\alpha\beta)\in H^2(M,K^2)$, multiplying by $s\in H^0(M,K^*)$ to get 
$$\sigma(\alpha,\beta)=s\tr(\alpha\beta)\in H^2(M,K)\cong \C.$$
The definition is very simple, the difficult part of \cite{Bot1} is proving the  vanishing of the Schouten bracket.

\begin{thm} \label{botta} Let $\sigma_+$ be the $I_+$ - Poisson structure defined by the generalized K\"ahler structure on ${\mathcal M}$. Then $\sigma_+/2$ is the canonical structure on the moduli space of $I_+$-stable bundles.
\end{thm}

\begin{prf} In the generalized K\"ahler setup, the Poisson structure $\sigma_+$ is defined by the $(0,2)$ part of $\omega_-$ under the antilinear identification $T^{1,0}\cong (\bar T^*)^{0,1}$ defined by the metric. 

A tangent vector to ${\mathcal M}$ is defined by $a\in \Omega^1(M,\lie{g})$ satisfying $d^+_Aa=0$, and this implies that $a^{0,1}\in \Omega^{0,1}(M,\End E)$ satisfies $\bar\partial_Aa^{0,1}=0\in \Omega^{0,2}(M,\End E)$, which is the tangent vector in the holomorphic setting -- it is a Dolbeault representative for a class in $H^1(M,\End E)$. 
The conjugate $a^{1,0}=\overline{a^{0,1}}$ defines a complex cotangent vector  by the linear form
$$b^{0,1}\mapsto \int_M\omega_+\wedge\tr (a^{1,0}\wedge b^{0,1})$$
and this is the  antilinear identification $T^{1,0}\cong (\bar T^*)^{0,1}$ on the moduli space. However $\omega_+\wedge a^{1,0}\in \Omega^{2,1}(M,\End E)$ is not a Dolbeault  representative for the Serre dual -- it is not $\bar\partial$-closed -- so to see concretely the canonical  Poisson structure we must find a good representative $(2,1)$ form.
\vskip .25cm
Now     from  $d_A^+a=0$ we have
$\omega_+\wedge d_A(a^{1,0}+a^{0,1})=0$ 
and from the horizontality condition $\omega_+\wedge d^c_+a=0$, we obtain 
$\omega_+\wedge d_A(a^{1,0}-a^{0,1})=0$
so putting them together
\begin{equation}
\omega_+\wedge\bar\partial_A a^{1,0}=0,\quad \omega_+\wedge\partial_A a^{0,1}=0
\label{infasd}
\end{equation}
\vskip .25cm
From Lemma \ref{hori} applied to $I_+$ and $I_-$ we have
$$\omega_{\pm} \wedge d_{\pm}^c a = \dstar a- d_{\pm}^c\omega_{\pm} \wedge a$$
so that since $d^c_-\omega_-=-d^c_+\omega_+$,
$$\omega_{-} \wedge d_{-}^c a = \omega_{+} \wedge d_{+}^c a +2 d_+^c\omega_{+} \wedge a.$$
If $a=a^{1,0}+d_A\psi$ is minus-horizontal then this equation tells us that 
$$0=\omega_{+} \wedge d_{+}^c d_A\psi+2d_+^c\omega_{+}\wedge (a^{1,0}+d_A\psi)$$
 since $a^{1,0}$ is plus-horizontal. We rewrite this as 
 \begin{equation}
2i\omega_{+}\wedge \bar\partial_A\partial_A \psi+2i\bar\partial\omega_{+}\wedge (a^{1,0}+\partial_A\psi)-2i\partial\omega_+\wedge\bar\partial_A\psi=0
\label{goodeq}
\end{equation}
using the fact that $\omega_+\wedge F=0$ where $F$ is the curvature of the connection $A$.
This gives, using $\bar\partial\partial\omega_+=0$ and (\ref{infasd}),
\begin{equation}
\bar\partial_A[\omega_+\wedge (a^{1,0}+\partial_A\psi)+\psi\partial\omega_+]=0
\label{Dolrep}
\end{equation}
Here, then, we have a 
 $\bar\partial$-closed form, and it  represents the dual of $[a^{0,1}]$ using the metric on ${\mathcal M}$ since, from Stokes' theorem,
$$\int_M[\omega_+\wedge \tr((a^{1,0}+\partial_A\psi)\wedge b^{0,1})+\partial\omega_+\wedge\tr(\psi b^{0,1})]= \int_M\omega_+\wedge \tr(a^{1,0}\wedge b^{0,1})-\int_M\omega_+\wedge\tr(\psi\partial_Ab^{0,1})$$
and the second term on the right hand side vanishes from (\ref{infasd}).
\vskip .25cm
Now  where the Poisson structure $s$ on $M$ is non-vanishing  we have a closed $2$-form $\beta_1-\bar\beta_2$ which   from Proposition \ref{4formulas} can be expressed as 
$2i\omega_+-(p-1)\bar\gamma/2+(p+1)\gamma/2.$ Since this is closed, and $\gamma$ is of type $(2,0)$,
$4id\omega_+=\partial p\wedge \bar\gamma-\bar\partial p\wedge \gamma$
and so 
\begin{equation}
4i\partial\omega_+=-\bar\partial p\wedge \gamma
\label{domega}
\end{equation}
 We can therefore rewrite the Dolbeault representative as 
$$\omega_+\wedge (a^{1,0}+\partial_A\psi)+i\psi\bar\partial p\wedge \gamma/4.$$
The canonical Poisson structure is therefore obtained by integrating over $M$ the form 
\begin{equation}
\tr[s(\omega_+\wedge (a_1^{1,0}+\partial_A\psi_1)+i\psi_1\bar\partial p\wedge \gamma/4)\wedge(\omega_+\wedge (a_2^{1,0}+\partial_A\psi_2)+i\psi_2\bar\partial p\wedge \gamma/4)]
\label{integrand}
\end{equation}
\vskip .25cm
Take the product of the  two expressions with an $\omega_+$ factor. 
For $(1,0)$ forms $a,b$, at each point
$[s(\omega_+\wedge a)]\wedge\omega_+\wedge b$ is a skew form on $T^{1,0}$ with values in $\Lambda^4T^*$ depending on a Hermitian form and a $(2,0)$ form $\gamma$ (recall from Proposition \ref{gprop} that $s\gamma=2i$). By $SU(2)$ invariance this must be a multiple of $\bar\gamma\wedge a\wedge b$ and a simple calculation shows that 
$$[s(\omega_+\wedge a)]\wedge\omega_+\wedge b=-i\frac{\omega_+^2}{\gamma\bar\gamma}\,\bar\gamma\wedge a\wedge b.$$
However from (\ref{eqA}) and $r=i(p^2-1)/4$ we see that
$$\frac{\omega_+^2}{\gamma\bar\gamma}=\frac{1}{8}(1-p^2).$$
But now from Proposition \ref{4formulas},  $\omega_-=p\omega_++i(p^2-1)\bar\gamma/4-i(p^2-1)\gamma/4$ and so
$$[s(\omega_+\wedge a)]\wedge\omega_+\wedge b=\frac{i}{2}\omega_-^{0,2}\wedge a\wedge b=\frac{i}{2}\omega_-\wedge a\wedge b$$
since $a$ and $b$ are of type $(1,0)$.
Thus the first two expressions contribute to the integral the term 
\begin{equation}
\frac{1}{2}\int_M\omega_- \wedge\tr(a_1^{1,0}+\partial_A\psi_1)\wedge (a_2^{1,0}+\partial_A\psi_2)
\label{two0}
\end{equation}
\vskip .25cm
The last two terms in (\ref{integrand}) give zero contribution because of the common $\bar\partial p$ factor. For the other terms, the relation $s\gamma=2i$ means that we are considering the integral of 
\begin{equation}
-\tr[\psi_1\bar\partial p\wedge\omega_+\wedge (a_2^{1,0}+\partial_A\psi_2)]/2+\tr[\psi_2\bar\partial p\wedge\omega_+\wedge (a_1^{1,0}+\partial_A\psi_1)]/2.
\label{crossterm}
\end{equation}
Take the first expression.
This no longer contains the singular term $\gamma$ so we can integrate over the manifold and  
using Stokes' theorem we get
\begin{equation}
\frac{1}{2}\int_M p\omega_+\wedge\tr(\bar\partial_A\psi_1\wedge (a_2^{1,0}+\partial_A\psi_2))+p\tr[\psi_1\bar\partial_A[\omega_+\wedge (a^{1,0}+\partial_A\psi)]
\label{integrate}
\end{equation}
Now from (\ref{Dolrep}) and (\ref{domega})
$$\bar\partial_A[\omega_+\wedge (a^{1,0}+\partial_A\psi)]=-\bar\partial_A(\psi\partial\omega_+)=\bar\partial_A\psi\wedge\bar\partial p\wedge \gamma/4i$$
 Using this we can write (\ref{integrate}) as 
\begin{equation}
 \frac{1}{2}\int_M p\omega_+\tr(\bar\partial_A\psi_1\wedge (a_2^{1,0}+\partial_A\psi_2))-\frac{i}{8}\int_M p\bar\partial p\wedge\tr(\bar\partial_A\psi_1\psi_2)\wedge \gamma
\label{integrate1}
\end{equation}
Now $\omega_-^{1,1}=p\omega_+$ and the first term integrates a $(1,1)$ form against $p\omega_+$ so we write this as 
\begin{equation}
\frac{1}{2}\int_M \omega_-\wedge\tr(\bar\partial_A\psi_1\wedge (a_2^{1,0}+\partial_A\psi_2))
\label{oneone}
\end{equation}
From Proposition \ref{4formulas}, we have $$\omega_-^{2,0}=-i(p^2-1)\gamma/4$$ so the last  term in (\ref{integrate1}) is
$$\frac{1}{4}\int_M \bar\partial \omega^{2,0}\wedge\tr(\bar\partial_A\psi_1\psi_2)$$
which using Stokes' theorem gives
$$\frac{1}{4}\int_M \omega_-^{2,0}\wedge\tr(\bar\partial_A\psi_1\bar\partial \psi_2)$$
which  we write as 
$$\frac{1}{4}\int_M \omega_-\wedge\tr(\bar\partial_A\psi_1\bar\partial \psi_2).$$
In the full integral there is another contribution of this form from the second term in (\ref{crossterm}) and adding all 
  terms in (\ref{integrand}) we obtain 
$$\frac{1}{2}\int_M\omega_-\wedge\tr(a_1^{1,0}+d_A\psi_1)\wedge\tr(a_2^{1,0}+d_A\psi_2).$$
Since $a_1^{1,0}+d_A\psi_1, a_2^{1,0}+d_A\psi_2$ are minus-horizontal representatives of $a_1^{1,0},a_2^{1,0}$ we see from the definition of $\tilde\omega_-$ that this is $\tilde\omega_-^{0,2}/2$ evaluated on those two vectors and hence  is half the Poisson structure defined by the bihermitian metric.
\end{prf}

\subsection{The generalized K\"ahler structure}

As we have seen, the bihermitian structure of $M^4$ naturally induces a similar structure  on the moduli space of instantons, but we only get a pair $J_1,J_2$ of commuting generalized complex structures by \emph{choosing} a $2$-form with $db=H$. In that respect $J_1,J_2$ are defined modulo a closed B-field but we can still extract some information about them. In particular the formula (\ref{J12}) shows that the real Poisson structures defined by $J_1$ and $J_2$, namely $\omega_+^{-1}\pm \omega_-^{-1}$, are unchanged by $b\mapsto b+B$. We shall determine the \emph{symplectic foliation} on ${\mathcal M}$ determined by these Poisson structures, which relates to the ``type" of the generalized complex structure as discussed by Gualtieri.

The symplectic foliation of a Poisson structure $
\pi$ is determined by  the subspace of the cotangent bundle annihilated by $\pi:T^*\rightarrow T$. From (\ref{J12}), in our case 
$$\ker \pi_1=\ker (I_++I_-),\quad \ker \pi_2=\ker (I_+-I_-)$$
where $I_+,I_-$ act on $T^*$.

Note that if $I_+a=I_-a$ then
$$[I_+,I_-]a=I_+I_-a-I_-I_+a=(I_+)^2a-(I_-)^2a=-a+a=0$$
so that $\ker (I_+-I_-)\subset \ker[I_+,I_-]$, and similarly if $I_+a=-I_-a$. It follows that if $I_+a=I_-a$, then $I_+(I_+a)=I_-(I_+a)$ since both sides are equal to $-a$. Thus $\ker \pi_1$ and $\ker \pi_2$ are complex subspaces of $\ker[I_+,I_-]$  (with respect to either structure). 

Now the kernel of $[I_+,I_-]$ is, from \ref{holp}, the kernel of the holomorphic Poisson structure $\sigma_+$ (or $\sigma_-$). But  Theorem \ref{botta} tells us that this is the canonical Poisson structure on ${\mathcal M}$. Its kernel is easily determined (see \cite{Bot1}). Recall that the  Poisson structure is defined, as a map from $(T^{1,0})^*$ to $T^{1,0}$, by the multiplication operation of the section $s$ of $K^*$:
$$s:H^1(M,\End E\otimes K)\rightarrow H^1(M,\End E).$$
If $D$ is the anticanonical divisor of $s$ then we have an exact sequence of sheaves
$$0\rightarrow{\mathcal O}_M(\End E\otimes K)\stackrel{s}\rightarrow {\mathcal O}_M(\End E)\rightarrow {\mathcal O}_D(\End E)\rightarrow 0$$
and the above is part of the long exact cohomology sequence. 
Since a stable bundle is simple, $H^0(M,\End E)$ is just the scalars, so the map
$H^0(M,\End E)\rightarrow H^0(D,\End E)$ just maps to the scalars. Hence the 
kernel of $\sigma_+$ is isomorphic from the exact sequence to $H^0(D,\End_0 E)$ under the connecting homomorphism: 
$$\delta_+: H^0(D,\End E)\rightarrow H^1(M,\End E\otimes K).$$
When $D$, an anticanonical divisor, is of multiplicity $1$ and smooth, it is an elliptic curve by the adjunction formula: $KD+D^2=2g-2$ implies $0=K(-K)+(-K)^2=2g-2$. Generically a holomorphic bundle on an elliptic curve is a sum of line bundles, and then the dimension of $H^0(D,\End_0 E)$ is $k-1$ if $\rk E=k$. Thus the real dimension of $\ker[I_+,I_-]$ is at least $2(k-1)$.
\vskip .25cm
Now the divisor $D$ is, by definition, the subset of $M$ on which $I_+=\pm I_-$, say $I_+=I_-$.
Thus the complex structure of the bundle $E$ determined by its ASD connection is the \emph{same} on $D$ for $I_+$ and $I_-$. So the same holomorphic section $u$ of $\End_0 E$ on $D$ maps  complex linearly in two different ways to the cotangent space of ${\mathcal M}$. To study these maps we should really say that there are real isomorphisms 
$$\alpha_{\pm}: H_{\pm}^1(M,\End E\otimes K)\rightarrow T^*_{[A]}$$
such that $\alpha_{\pm}$ is $I_{\pm}$-complex linear.

\begin{prp} \label{residue} $\alpha_+\delta_+=\alpha_-\delta_-$
\end{prp}

\begin{prf} Recall how the connecting homomorphism is defined in Dolbeault terms, for the moment in the case where $D$ has multiplicity one: we have a holomorphic section $u$ of $\End_0 E$ on $D$, and then extend using a partition of unity to a $C^{\infty}$ section $\tilde u$ on $M$. Then since $u$ is holomorphic on $D$, $\bar\partial \tilde u$ is divisible by $s$, the section of $K^*$ whose divisor is $D$. Then
$\delta(u)$ is represented by the $(2,1)$-form $s^{-1}\bar\partial_A\tilde u$. 

Let $a\in T_{[A]}$ be a tangent vector to the moduli space, so $a\in\Omega^1(M,\End E)$ and satisfies $d_A^+a=0$. So $\bar\partial_A a^{0,1}=0$ and we evaluate the cotangent vector $\delta_+(u)$ on $a$ to get
$$\int_M\tr(s^{-1}\bar\partial_A\tilde u \wedge a).$$
But $s^{-1}=\gamma/2i$ so this is
 $$\frac{1}{2i}\int_M\gamma\tr(\bar\partial_A\tilde u\wedge  a).$$
Away from the divisor $D$, we have
$$\bar\partial(\gamma\wedge\tr(\tilde u a))=\gamma\wedge \tr(\bar\partial_A\tilde u \wedge a)$$
since both $\gamma$ and $a$ are $\bar\partial$-closed. By Stokes' theorem the integral is reduced to an integral around the unit circle bundle of the normal bundle of $D$ and from there to an integral over $D$. In fact, if $\gamma$ has a simple pole along $D$ then locally
$$\gamma=f(z_1,z_2)\frac{dz_1\wedge dz_2}{z_1}$$
where $z_1=0$ is the equation of $D$. The holomorphic one-form $f(0,z_2)dz_2$ is then globally defined on $D$ -- the \emph{residue} $\gamma_0$ of the meromorphic $2$-form. This residue is the same for $I_+$ and $I_-$ (from Proposition \ref{4formulas} the meromorphic form for $I_-$ is $-2i\omega_+-(p-1)/2\gamma+(p+1)/2\bar\gamma$ and $p=-1$ on $D$). 
 Thus the integral becomes
$$\frac{1}{2i}\int_D\gamma_0\wedge \tr(ua).$$
This is defined entirely in terms of the data on $D$ and so is the same for $I_+$ and $I_-$.
\vskip .25cm
When the divisor has multiplicity $d$, the section $u$  extends holomorphically to the $(d-1)$-fold formal neighbourhood of the curve and our $C^{\infty}$ extension  must agree with this. The result remains true. (Note that the discussion  of Poisson surfaces and moduli spaces via the residue is the point of view advanced in Khesin's work \cite{Kh}.)
\end{prf}
\vskip .25cm
\begin{cor} The two real Poisson structures $\pi_1,\pi_2$ defined by the generalized complex structures $J_1,J_2$  on the moduli space ${\mathcal M}$ of $SU(k)$ instantons  have kernels of dimension $0$ and $\ge 2(k-1)$.
\end{cor}
\begin{prf}
We saw at the beginning of the Section that if $I_+a=I_-a$ then $[I_+,I_-]a=0$. Proposition \ref{residue} shows that $I_+$ and $I_-$ agree on the kernel of $[I_+,I_-]$, so that $\ker(I_+-I_-)=\ker[I_+,I_-]$. 

Now $\ker(I_+-I_-)$ is the kernel of Poisson structure $\pi_1$ say, which is isomorphic to $H^0(D,\End_0 E)$ and has, as we have seen, at least $2(k-1)$ real dimensions. The other Poisson structure $\pi_2$ has kernel $\ker(I_++I_-)$. But this also lies in the kernel of $[I_+,I_-]$ so $I_+a=I_-a$. With $I_+a=-I_-a$ this means $a=0$.
\end{prf}
\vskip .25cm
The generalized complex structure $J_2$ on ${\mathcal M}$ where the kernel of the Poisson structure is zero is therefore of the form $\exp(B+i\omega)$ and it is tempting to associate it to the generalized complex structure  of symplectic type on $M^4$. However, as we have seen, there appears to be no way to naturally associate or even define these   structures, since the $2$-form $b$ does not descend to the moduli space.

\subsection{Examples of symplectic leaves}

We saw in the previous section that the symplectic leaves of $\pi_1$ are the same as the the symplectic leaves of the canonical complex Poisson structure on ${\mathcal M}$.

The simplest example is to take $\CP^2$ with the anticanonical divisor defined by a triple line : $D=3L$. The moduli space of stable rank $2$ bundles with $c_2=2$ has dimension $4\times 2-3=5$ and has a very concrete description. Such a bundle $E$ is trivial on a general projective line but jumps to ${\mathcal O}(1)\oplus {\mathcal O}(-1)$ on the lines which are tangent to a nonsingular conic $C_E$. The moduli space ${\mathcal M}$ is then just the space of non-singular conics, which is a homogeneous space of $PGL(3,\C)$.

The subgroup preserving $L$ (the line at infinity say) is the affine group $A(2)$ and if it preserves the Poisson structure it fixes $dz_1\wedge dz_2$.  Hence the $5$-dimensional unimodular affine group $SA(2)$ acts on ${\mathcal M}$ preserving the Poisson structure. The subgroup $G$ which fixes the conic $z_1z_2=a$ consists of the transformations $(z_1,z_2)\mapsto (\lambda z_1,\lambda^{-1}z_2)$ so for each $a$, the orbit of the conic under $SA(2)$ is isomorphic to the $4$-dimensional quotient $SA(2)/G$.  These orbits are the generic symplectic leaves of the Poisson structure, and thus are homogeneous symplectic and hence isomorphic to coadjoint orbits. 

In fact if $z\mapsto Az+b$ is in the Lie algebra of $SA(2)$   then $G$ is the stabilizer of the linear map $f(A,b)=A_{11}$ so that $SA(2)/G$ is the orbit of $f$ in the dual of the Lie algebra.
 This deals with conics which meet $L$ in two points. The ones which are tangential to $L$ (i.e. the bundles for which $L$ is a jumping line) are parabolas: e.g. $z_1^2=z_2$. The identity component of the stabilizer of this is the one-dimensional group 
$(z_1,z_2)\mapsto (z_1+c,2cz_1+z_2+c^2)$ and this stabilizes the linear map
$(A,b)\mapsto A_{21}+4b_1$, so we again have a coadjoint orbit.

\vskip .25cm
In general, the symplectic leaves are roughly given by the bundles $E$ on $M$ which restrict to the same bundle on the anticanonical divisor $D$. ``Roughly", because we are looking at equivalence classes and a stable bundle on $M$ may not restrict to a stable bundle on $D$, so there may not be a well-defined map from ${\mathcal M}$ to a Hausdorff moduli space. On the other hand this is the quotient space of a (singular) foliation  so we don't expect that.

When $D$ is the triple line $3L$ in $\CP^2$ there is an alternative way of describing these leaves. On a generic line $E$ is trivial and the sections along that line define the fibre of a vector bundle $F$ on the dual plane, outside the curve $J$ of jumping lines. If we take a section of $E$ on $L$ we can try and extend it to the first order neighbourhood of $L$. Since the normal bundle to $L$ is ${\mathcal O}(1)$ there is an exact sequence of sheaves for sections on the $n$-th order neighbourhood:
$$0\rightarrow {\mathcal O}(E(-n))\rightarrow {\mathcal O}^{(n)}(E)\rightarrow {\mathcal O}^{(n-1)}(E)\rightarrow 0.$$
Since $H^0(\CP^1, {\mathcal O}(-1))=H^1(\CP^1, {\mathcal O}(-1))=0$, any section has a unique extension to the first order neighbourhood: this defines a \emph{connection} on $F$. The extension to the second order neighbourhood is obstructed since $H^1(\CP^1, {\mathcal O}(-2))\cong \C$ and this obstruction is the \emph{curvature} of the connection (see \cite{Hurt} for details of this twistorial construction).

What it means is that if $L$ is not a jumping line, then $E$ restricted to $3L$ is essentially the curvature of the connection on $F$ at the point $\ell$ in the dual plane defined by the line $L$, and the symplectic leaves are obtained by fixing the equivalence class of the curvature at that point. The curvature acquires a double pole on $J$.

 From this point of view, the case $k=2,c_1=0,c_2=2$ concerns an $SO(3,\C)$-invariant connection on a rank $2$-bundle on the complement of a conic, and this is essentially the Levi-Civita connection of $\RP^2$ complexified. This is an $O(2)$-connection which becomes an $SO(2)$ connection on $S^2$ with curvature
$$\frac{dz\wedge d\bar z}{(1+\vert z\vert^2)^2}.$$
So the bundle on $D$ is equivalent to the transform of the complexification of this by a projective transformation. If the dual conic is defined by the symmetric $3\times 3$ matrix $Q_{ij}$ and $x$ is a vector representing $\ell$ then the curvature is
$$\frac{(\det Q)^{2/3}}{Q(x,x)^2}.$$
The symplectic leaves are then given by the equation $\det Q=aQ(x,x)^3$ for varying $a$.

\section{A quotient construction}

It is well-known that the moduli space of instantons on a hyperk\"ahler $4$-manifold is hyperk\"ahler and this can be viewed as an example in infinite dimensions of a hyperk\"ahler quotient -- the quotient of the space of all connections by the action of the group of gauge transformations. One may ask if, instead of the painful integration by parts that we did in the previous sections, there is a cleaner way of viewing the definition of a generalized K\"ahler structure on ${\mathcal M}$. The problem is that such a quotient would have to encompass not only the hyperk\"ahler quotient but also the ordinary K\"ahler quotient, and in finite dimensions these are very different -- the dimension of the quotient in particular is different! 

We offer next an example of a generalized K\"ahler quotient which could be adapted to replace the  differential geometric arguments in the previous sections for the case of a torus or K3, and at least gives another reason why the calculations should hold.  It also brings out in a natural way the frustrating feature that the $2$-form $b$ does not descend in general to the quotient.
\vskip .25cm

We suppose the  generalized K\"ahler structure  is even and is given by global forms 
$\rho_1=\exp {\beta_1},  \rho_2=\exp{\beta_2}$ 
where $\beta_1,\beta_2$ are closed complex forms on a real manifold $M$ of dimension $4k$. This is the test situation we have been considering throughout this paper. From Lemma \ref{commute}, the compatibility ($J_1J_2=J_2J_1$) is equivalent to 
$$(\beta_1-\beta_2)^{k+1}=0,\qquad (\beta_1-\bar\beta_2)^{k+1}=0.$$

Now suppose a Lie group $G$ acts, preserving the forms $\beta_1,\beta_2$, and giving complex moment maps $\mu_1,\mu_2$.  To make a quotient, we would like to take the joint zero set of $\mu_1$ and $\mu_2$ and divide by the group $G$, but these are two \emph{complex} functions so if they were generic we would get as a quotient a manifold of dimension $\dim M-5\dim G$ instead of  $\dim M-4\dim G$.

To avoid this, we need to assume that $\beta_1,\beta_2,\bar\beta_1,\bar\beta_2$ are linearly dependent over $\R$. 

\begin{rmk} If we were trying to set up the moduli space  of instantons as a quotient of the space of all connections on a K3 surface or a torus, the following lemma links the condition of linear dependence of the moment maps to the necessity to choose a self-dual $b$. 

\begin{lem} If $\dim M=4$, then $\beta_1,\beta_2,\bar\beta_1,\bar\beta_2$ are linearly independent over $\R$ at each point  if and only if $b$ is self-dual.
\end{lem}

\begin{lemprf} From Proposition \ref{4formulas} we have
\begin{eqnarray*}
\beta_1+\bar\beta_1&=&2b-(p-1)(\gamma+\bar\gamma)/2\\
\beta_2+\bar\beta_2&=&2b-(p+1)(\gamma+\bar\gamma)/2\\
-i(\beta_1-\bar\beta_1)&=&2\omega_+-(p-1)i(\gamma-\bar\gamma)/2\\
-i(\beta_2-\bar\beta_2)&=&2\omega_+-(p+1)i(\gamma-\bar\gamma)/2\\
\end{eqnarray*}
We can easily solve these for $b,\omega_+,\gamma+\bar\gamma,i(\gamma-\bar\gamma)$ in terms of the $\beta_i$. 

If $b$ is self-dual, it is a real linear combination of $\omega_+,\gamma+\bar\gamma,i(\gamma-\bar\gamma)$ since $\gamma$ is of type $(2,0)$ relative to $I_+$, hence we get a linear relation amongst the left hand sides.

Conversely, a linear relation among the left hand sides  will express $b$ in terms of $\omega_+,\gamma+\bar\gamma,i(\gamma-\bar\gamma)$ unless it is of the form
$$(\beta_1+\bar\beta_1)-(\beta_2+\bar\beta_2)+i\lambda(\beta_1-\bar\beta_1)+i\mu (\beta_2-\bar\beta_2)=0.$$
But the $(1,1)$ component of this is $2(\mu-\lambda)\omega_+$ so $\lambda=\mu$ and then the relation can be written
$$(1+i\lambda)(\beta_1-\bar\beta_2)+(1-i\lambda)(\bar\beta_1-\beta_2)=0$$
but $(\bar\beta_1-\beta_2)$ is of type $(2,0)$ relative to $I_-$ so this is impossible.
\end{lemprf}

We see that the condition for $b$ to define $\tilde b$ on the moduli space ${\mathcal M}$ is related to the  linear dependence issue of the moment maps.
\end{rmk}

Returning to the general case, for  each vector field $X$ from the Lie algebra of $G$ we have $i_X\beta_i=d\mu_i$ and so $\beta_i$ restricted to $\mu_1=0=\mu_2$ is annihilated by $X$, and invariant under the group and hence is the pullback of a form $\tilde\beta_i$ on the quotient, which is also closed. 

In the bihermitian interpretation, $\bar\beta_1-\bar\beta_2$ is a non-degenerate $(2,0)$-form relative to $I_+$ -- a holomorphic symplectic form -- and the quotient can then be identified  with the holomorphic symplectic quotient. In particular if the complex dimension of the quotient is $2m$ then $(\tilde\beta_1-\tilde\beta_2)^{m+1}=0$ and $(\tilde\beta_1-\tilde\beta_2)^{m}\ne 0$. Similarly $(\bar\beta_1-\beta_2)$ is $(2,0)$ with respect to $I_-$ and we get the same property for $(\tilde\beta_1-\bar{\tilde \beta_2})$. From Lemma \ref{commute} we have a generalized K\"ahler structure on the quotient.
\vskip .25cm
Note that  in this generic case the Poisson structures on the quotient are non-degenerate.


\begin{thebibliography}{11}
 %
 \bibitem{AB}
 M. Abouzaid and M. Boyarchenko, {\it Local structure of generalized complex manifolds},  {\bf math.DG/0412084}
%
\bibitem{AGG}
V. Apostolov, P. Gauduchon and G. Grantcharov, {\it Bihermitian structures on complex surfaces}, Proc. London Math. Soc. {\bf 79} (1999), 414--428
%
\bibitem{Bart}
C. Bartocci and E. Macr\`\i , {\it Classification of Poisson surfaces},  Communications in Contemporary Mathematics (to appear) {\bf math.AG/0402338}.
%
\bibitem{Bot1}
F. Bottacin,
{\it Poisson structures on moduli spaces of sheaves over Poisson surfaces},
Invent. Math. {\bf 121} (1995), 421--436.
%
\bibitem{Buch}
N. P. Buchdahl, {\it Hermitian-Einstein connections and stable vector bundles over compact complex surfaces}, Math. Ann. {\bf 280} (1988) 625--648.
%
\bibitem{R}
S. J. Gates, C. M. Hull and M. Ro\v cek, {\it Twisted multiplets and new supersymmetric nonlinear $\sigma$-models}.  Nuclear Phys. B  {\bf 248}  (1984),   157--186.
%
 \bibitem{Gu}
 M. Gualtieri, {\it Generalized complex geometry}, 
{\bf math.DG/0401221}
%
 \bibitem{Hit}
 N. J. Hitchin, {\it Generalized Calabi-Yau manifolds}, Q. J. Math. {\bf 54} (2003) 281--308
%
\bibitem{Hurt}
J. Hurtubise, {\it Twistors and the geometry of bundles over $P_2(C)$}, Proc. London Math. Soc. {\bf 55} (1987) 450--464.
%
\bibitem{Kh} 
 B. Khesin and A. Rosly,
 {\it Symplectic geometry on moduli spaces of holomorphic bundles over complex surfaces},  
 The Arnoldfest (Toronto, ON, 1997),  311--323, 
 Fields Inst. Commun., {\bf 24} 
 Amer. Math. Soc., Providence, RI (1999). 

\bibitem{LY}
J. Li and S-T. Yau, {\it Hermitian Yang-Mills connections on non-K\"ahler manifolds}, in ``Mathematical aspects of string theory", World Scientific (1987)
%
\bibitem{LT}
 M. L\"ubke and A. Teleman, ``The Kobayashi-Hitchin correspondence", World Scientific, Singapore (1995) 
\end{thebibliography}
\end{document}